\nonstopmode
\documentclass[a4paper, 10pt]{amsart}
\usepackage{latexsym}
\usepackage{fancyhdr}
\usepackage{longtable}
\usepackage{amsmath, amssymb}
\usepackage[ansinew]{inputenc}
\usepackage[all]{xy}
\usepackage{verbatim}
\usepackage{psfrag}
\usepackage{epsfig}
\swapnumbers
\theoremstyle{plain}
\newtheorem*{lemma*}{Lemma}
\newtheorem{lemma}[subsection]{Lemma}
\newtheorem*{theorem*}{Theorem}
\newtheorem{theorem}[subsection]{Theorem}
\newtheorem*{proposition*}{Proposition}
\newtheorem{proposition}[subsection]{Proposition}
\newtheorem*{corollary*}{Corollary}
\newtheorem{corollary}[subsection]{Corollary}
\theoremstyle{definition}
\newtheorem*{definition*}{Definition}

\newtheorem*{example*}{Example}

\newtheorem*{remark*}{Remark}
\newtheorem*{remarks*}{Remarks}

\numberwithin{equation}{subsection}
\newenvironment{demo}[1]{\par\smallskip\noindent{\bf #1.}}{\par\smallskip}

\newcommand{\East}[2]{-\raisebox{0.1pt}{$\mkern-16mu\frac{\;\;#1\;}{\;\;#2\;}\mkern-16mu$}\to}

\let\on=\operatorname

\title[Denjoy--Carleman mappings]
{The Convenient Setting for non-Quasianalytic Denjoy--Carleman
        Differentiable Mappings}

\author
[A.~Kriegl, P.W.~Michor, A.~Rainer]
{Andreas Kriegl, Peter W. Michor, and Armin Rainer}

\address{Andreas Kriegl: Fakult\"at f\"ur Mathematik, Universit\"at Wien, 
        Nordbergstrasse~15, A-1090 Wien, Austria}
\email{andreas.kriegl@univie.ac.at}

\address{Peter W. Michor: Fakult\"at f\"ur Mathematik, Universit\"at Wien, 
        Nordbergstrasse~15, A-1090 Wien, Austria}
\email{peter.michor@univie.ac.at}

\address{Armin Rainer: Fakult\"at f\"ur Mathematik, Universit\"at Wien, 
        Nordbergstrasse~15, A-1090 Wien, Austria}
\email{armin.rainer@univie.ac.at}

\thanks{PM was supported by FWF-Project P~21030-N13.
        AR was supported by FWF-Projects P19392 \& J2771'}

\begin{document}
        
        \begin{abstract}
                For Denjoy--Carleman differentiable function classes $C^M$ 
                where the weight sequence 
                $M=(M_k)$ is logarithmically convex, stable under derivations, 
                and non-quasianalytic of moderate growth, we prove the following: 
                A mapping is $C^M$ if it maps $C^M$-curves to $C^M$-curves.
                The category of $C^M$-mappings is cartesian closed in the sense that 
                $C^M(E,C^M(F,G))\cong C^M(E\times F, G)$ for convenient vector spaces. 
                Applications to manifolds of mappings are given: The group of
                $C^M$-diffeomorphisms is a $C^M$-Lie group but not better. 
        \end{abstract}

        \maketitle
        
        \section{\label{nmb:1} Introduction}
        Denjoy--Carleman differentiable functions form spaces of functions between
        real analytic and $C^\infty$. They are described by growth conditions on
        the Taylor expansions, see \thetag{\ref{nmb:2.1}}. 
        Under appropriate conditions the
        fundamental results of calculus still hold:
        Stability under differentiation, composition, solving ODEs, applying the
        implicit function theorem. 
        See Section \thetag{\ref{nmb:2}} for a review of Denjoy--Carleman differentiable
        functions, which is summarized in Table 1. 
        
        In \cite{Kriegl82}, \cite{Kriegl83}, \cite{FK88}, \cite{KrieglNel85}, 
        \cite{KrieglMichor90}, see \cite{KM97} for a comprehensive presentation,
        convenient calculus was developed for $C^\infty$,
        holomorphic, and real analytic functions: see appendix \thetag{\ref{nmb:7}}, \thetag{\ref{nmb:8}},
        \thetag{\ref{nmb:9}} for a short overview of the essential results. 
        
        In this paper we develop the convenient calculus for Denjoy--Carleman
        classes $C^M$ where the weight sequence 
        $M=(M_k)$ is logarithmically convex, stable under derivations, and 
        non-quasianalytic of moderate growth 
        (this holds for all Gevrey differentiable functions
        $G^{1+\delta}$ for $\delta>0$).
        By `convenient calculus' we mean that the following
        theorems are proved: 
        A mapping is $C^M$ if it maps $C^M$-curves to $C^M$-curves, see \thetag{\ref{nmb:3.9}}; 
        this is wrong in the quasianalytic case, see \ref{nmb:3.12}. 
        The category of $C^M$-mappings is cartesian closed in the sense that 
        $C^M(E,C^M(F,G))\cong C^M(E\times F, G)$ for convenient vector spaces, see
        \thetag{\ref{nmb:5.3}}; this is wrong for weight sequences of non-moderate growth, 
        see \thetag{\ref{nmb:5.4}}. The uniform boundedness principle holds for linear mappings
        into spaces of $C^M$-mappings. 
        
        For the quasianalytic case we hope for results similar to the real analytic
        case, but the methods have to be different. This will be taken
        up in another paper. 
        
        In chapter \thetag{\ref{nmb:6}} some applications to manifolds of mappings are given: The group of
        $C^M$-diffeomorphisms is a $C^M$-Lie group but not better.

        \section{\label{nmb:2} Review of Denjoy--Carleman differentiable functions}
        
        \subsection{\label{nmb:2.1}Denjoy--Carleman classes $C^M(\mathbb R^n,\mathbb R)$ 
                of differentiable functions}
        
        We mainly follow \cite{Thilliez08} (see also the references therein).
        We use $\mathbb{N} = \mathbb{N}_{>0} \cup \{0\}$.
        For each multi-index $\alpha=(\alpha_1,\ldots,\alpha_n) \in \mathbb{N}^n$, we write
        $\alpha!=\alpha_1! \cdots \alpha_n!$, $|\alpha|= \alpha_1 +\cdots+ \alpha_n$, and 
        $\partial^\alpha=\partial^{|\alpha|}/\partial x_1^{\alpha_1} \cdots \partial x_n^{\alpha_n}$.

        Let $M=(M_k)_{k \in \mathbb{N}}$ be an increasing sequence ($M_{k+1}\ge M_k$) 
        of positive real numbers with $M_0=1$.
        Let $U \subseteq \mathbb{R}^n$ be open.
        We denote by $C^M(U)$ the set of all $f \in C^\infty(U)$ such that, for all
        compact $K \subseteq U$,
        there exist positive constants $C$ and $\rho$ such that
        \begin{equation}\label{CM}
                |\partial^\alpha f(x)| \le C \, \rho^{|\alpha|} \, |\alpha|! \, M_{|\alpha|}
        \end{equation}
        for all $\alpha \in \mathbb{N}^n$ and $x \in K$.
        The set $C^M(U)$ is a \emph{Denjoy--Carleman class} of functions on $U$.
        If $M_k=1$, for all $k$, then $C^M(U)$ coincides with the ring $C^\omega(U)$
        of real analytic functions
        on $U$. In general, $C^\omega(U) \subseteq C^M(U) \subseteq C^\infty(U)$.
        
        We assume that $M=(M_k)$ is \emph{logarithmically convex}, i.e.,
        \begin{equation} \label{logconvex}
                M_k^2 \le M_{k-1} \, M_{k+1} \quad \text{ for all } k,
        \end{equation}
        or, equivalently, $M_{k+1}/M_k$ is increasing.
        Considering $M_0=1$, we obtain that also $(M_k)^{1/k}$ is increasing and
        \begin{equation}
                M_l \, M_k\le M_{l+k} \quad \text{ for all }l,k\in \mathbb{N}.
                \label{logconvex1}
        \end{equation}
        We also get (see \thetag{\ref{nmb:2.9}})
        \begin{equation}
                M_1^k \, M_k\ge M_j\, M_{\alpha_1} \cdots M_{\alpha_j} \quad \text{ for all }
                \alpha_i\in \mathbb{N}_{>0}, \alpha_1+\dots+\alpha_j = k.
                \label{logconvex2}
        \end{equation}
        Let $M=(M_k)$ be logarithmically convex. 
        Then $M_k'= M_k/M_0\,M_1^k\ge 1$ is increasing by \eqref{logconvex2}, 
        logarithmically convex, and 
        $C^M(U)=C^{M'}(U)$ for all $U$ open in $\mathbb R^n$ by \eqref{incl}. 
        So without loss we assumed at the beginning that $M$ is increasing.
        
        Hypothesis \eqref{logconvex} implies that $C^M(U)$ is a ring, for all open
        subsets $U \subseteq \mathbb{R}^n$, 
        which can easily be derived from \eqref{logconvex1} by means of Leibniz's rule.
        Note that definition \eqref{CM} makes sense also for mappings 
        $U\to \mathbb R^p$.
        For $C^M$-mappings, \eqref{logconvex} guarantees stability under composition
        (\cite{Roumieu62/63}, see also \cite[4.7]{BM04}; 
        a proof is also contained in the end of the
        proof of \thetag{\ref{nmb:3.9}}).

        A further consequence of \eqref{logconvex} is the inverse function theorem
        for $C^M$ (\cite{Komatsu79}; for a proof see also \cite[4.10]{BM04}): Let
        $f : U \to V$ be a $C^M$-mapping between open subsets $U,V \subseteq \mathbb{R}^n$.
        Let $x_0 \in U$.
        Suppose that the Jacobian matrix $(\partial f/\partial x)(x_0)$ is invertible. Then
        there are neighborhoods $U'$ of $x_0$, $V'$ of
        $y_0 := f(x_0)$ such that $f:U'\to V'$ is a $C^M$-diffeomorphism.

        Moreover, \eqref{logconvex} implies that $C^M$ is closed under solving ODEs
        (due to \cite{Komatsu80}):
        Consider the initial value problem
        \[
                \frac{d x}{d t} = f(t,x), \quad x(0) = y,
        \]
        where $f : (-T,T) \times \Omega \to \mathbb{R}^n$, $T>0$, and $\Omega \subseteq \mathbb{R}^n$ is
        open.
        Assume that $f(t,x)$ is Lipschitz in $x$, locally uniformly in $t$. 
        Then for each relative compact open
        subset $\Omega_1 \subseteq \Omega$ there exists $0 < T_1 \le T$ such that for
        each $y \in \Omega_1$ there is a unique solution $x=x(t,y)$ on the
        interval $(-T_1,T_1)$.
        If $f : (-T,T) \times \Omega \to \mathbb{R}^n$ is a $C^M$-mapping then the solution 
        $x : (-T_1,T_1) \times \Omega_1 \to \mathbb{R}^n$
        is a $C^M$-mapping as well.

        Suppose that $M=(M_k)$ and $N=(N_k)$ satisfy $M_k \le C^k \, N_k$, for all $k$
        and a constant $C$, or equivalently,
        \begin{equation} \label{incl}
                \sup_{k \in \mathbb{N}_{>0}} \Big(\frac{M_k}{N_k}\Big)^{\frac{1}{k}} < \infty.
        \end{equation}
        Then, evidently $C^M(U) \subseteq C^N(U)$. The converse
        is true as well (if \eqref{logconvex} is assumed): 
        One can prove that there exists $f \in C^M(\mathbb{R})$ such that 
        $|f^{(k)}(0)| \ge k! \, M_k$ for all $k$ (see \cite[Theorem 1]{Thilliez08}).
        So the inclusion $C^M(U) \subseteq C^N(U)$ implies \eqref{incl}. 
        
        Setting $N_k=1$ in \eqref{incl} yields that $C^\omega(U) = C^M(U)$ if and only if 
        \[
                \sup_{k \in \mathbb{N}_{>0}} (M_k)^{\frac{1}{k}} < \infty.
        \] 
        Since $(M_k)^{1/k}$ is increasing (by logarithmic convexity), 
        the strict inclusion $C^\omega(U) \subsetneq C^M(U)$ is equivalent to 
        \[
                \lim_{k \to \infty} (M_k)^{\frac{1}{k}} = \infty.
        \] 
        
        We shall also assume that $C^M$ is stable under derivation, which is
        equivalent to the following condition
        \begin{equation} \label{der}
                \sup_{k \in \mathbb{N}_{>0}} \Big(\frac{M_{k+1}}{M_k}\Big)^{\frac{1}{k}} < \infty.
        \end{equation}
        Note that the first order partial derivatives of elements in $C^M(U)$
        belong to $C^{M^{+1}}(U)$,
        where $M^{+1}$ denotes the shifted sequence 
        $M^{+1} = (M_{k+1})_{k \in \mathbb{N}}$.
        So the equivalence follows from \eqref{incl}, by replacing $M$ with
        $M^{+1}$ and $N$ with $M$.
        
        \subsection*{Definition} 
        By a {\it DC-weight sequence} we mean a sequence 
        $M=(M_k)_{k \in \mathbb{N}}$ of positive numbers with $M_0=1$ 
        which is monotone increasing
        ($M_{k+1}\ge M_k$), logarithmically convex
        \eqref{logconvex}, and satisfies \eqref{der}.
        Then $C^M(U,\mathbb R)$ is a differential ring, and the class of
        $C^M$-functions is stable under compositions.
        DC stands for Denjoy-Carleman and also for derivation closed.
        
        \subsection{\label{nmb:2.2} Quasianalytic function classes}
        Let $\mathcal{F}_n$ denote the ring of formal power series in $n$ variables (with
        real or complex coefficients).
        For a sequence $M_0=1, M_1,M_2,\dots>0$, 
        we denote by $\mathcal{F}_n^M$ the set of elements 
        $F=\sum_{\alpha \in \mathbb{N}^n} F_\alpha \, x^\alpha$ of $\mathcal{F}_n$ for which there 
        exist positive constants $C$ and $\rho$ such that
        \[
                |F_\alpha| \le C \, \rho^{|\alpha|} \, M_{|\alpha|}
        \]
        for all $\alpha \in \mathbb{N}^n$.
        A class $C^M$ is called \emph{quasianalytic} if, for open connected $U
        \subseteq \mathbb{R}^n$ and all $a \in U$, the Taylor series homomorphism
        \[
                T_a : C^M(U) \to \mathcal{F}_n^M,
                ~f \mapsto T_a f(x) = \sum_{\alpha \in \mathbb{N}^n} \frac{1}{\alpha!} \, \partial^\alpha f(a) \, x^\alpha
        \]
        is injective.
        By the Denjoy--Carleman theorem (\cite{Denjoy21}, \cite{Carleman26}),
        {\it the following statements are equivalent:
                \begin{enumerate}
                \item $C^M$ is quasianalytic.
                \item $\sum_{k=1}^\infty \frac1{m_k} = \infty$
                        where $m_k= \inf\{(j!\,M_j)^{1/j}: j\ge k\}$ is the increasing minorant of 
                        $(k!\,M_k)^{1/k}$.
                \item $\sum_{k=1}^\infty (\frac1{M_k^*})^{1/k} = \infty$
                        where $M^*_k= \inf\{(j!\,M_j)^{(l-k)/(l-j)}(l!\,M_l)^{(k-j)/(l-j)}: j\le k\le l, j<l\}$ 
                        is the logarithmically convex minorant of $k!\,M_k$.
                \item 
                        $\sum_{k=0}^\infty \frac{M^*_k}{M^*_{k+1}}=\infty$.
                \end{enumerate}
        }
        For contemporary proofs see for instance \cite[1.3.8]{Hoermander83I} 
        or \cite[19.11]{Rudin87}. 
        
        Suppose that $C^\omega(U) \subsetneq C^M(U)$ and $C^M(U)$ 
        is quasianalytic and logarithmically convex. 
        Then $T_a : C^M(U) \to \mathcal{F}_n^M$ is not surjective. 
        This is due to Carleman \cite{Carleman26}; 
        an elementary proof can be found in \cite[Theorem 3]{Thilliez08}.
        
        \subsection{\label{nmb:2.3} Non-quasianalytic function classes}
        If $M$ is a DC-weight sequence which is not quasianalytic, then there are 
        $C^M$ partitions of unity. Namely, there exists a $C^M$ function $f$ on $\mathbb R$ 
        which does not vanish in any neighborhood of 0 but which has vanishing
        Taylor series at 0. Let $g(t)=0$ for $t\le 0$ and $g(t)=f(t)$ for $t>0$.
        From $g$ we can construct $C^M$ bump functions as usual. 
        
        \subsection{\label{nmb:2.4} Strong non-quasianalytic function classes}
        Let $M$ be a DC-weight sequence with $C^\omega(U, \mathbb R) \subsetneq C^M(U,\mathbb R)$. 
        Then the mapping $T_a : C^M(U,\mathbb R) \to \mathcal{F}_n^M$
        is surjective, for all $a \in U$,
        if and only if there is a constant $C$ such that
        \begin{equation} \label{strnonqa}
                \sum_{k=j}^\infty \frac{M_k}{(k+1) \, M_{k+1}} \le C \frac{M_j}{M_{j+1}} \quad
                \text{for any integer}~ j \ge 0.
        \end{equation}
        See \cite{Petzsche88} and references therein.
        \eqref{strnonqa} is called \emph{strong non-quasianalyticity} condition.
        
        \subsection{\label{nmb:2.5} Moderate growth}
        A DC-weight sequence $M$ has
        \emph{moderate growth} if
        \begin{equation} \label{mgrowth}
                \sup_{j,k \in \mathbb{N}_{>0}} \Big(\frac{M_{j+k}}{M_j \, M_k}\Big)^{\frac{1}{j+k}} <
                \infty.
        \end{equation}
        Moderate growth implies derivation closed. 
        
        Moderate growth together with 
        strong non-quasianalyticity \eqref{strnonqa} is called
        \emph{strong regularity}: Then
        a version of Whitney's
        extension theorem holds for
        the corresponding function classes (e.g.\ \cite{BBMT91}).
        
        \subsection{\label{nmb:2.6} Gevrey functions}
        Let $\delta >0$ and put $M_k=(k!)^\delta$, for $k \in \mathbb{N}$.
        Then $M=(M_k)$ is strongly regular. The corresponding class $C^M$ of functions is the
        \emph{Gevrey class} $G^{1+\delta}$.
        
        \subsection{\label{nmb:2.7}More examples}
        Let $\delta >0$ and put $M_k=(\log(k+e))^{\delta \, k}$, for $k \in \mathbb{N}$.
        Then $M=(M_k)$ is quasianalytic for $0 < \delta \le 1$ and non-quasianalytic 
        (but not strongly) for $\delta > 1$. In any case $M$ is of moderate growth. 
        
        Let $q > 1$ and put $M_k= q^{k^2}$, for $k \in \mathbb{N}$. 
        The corresponding $C^M$-functions are called \emph{$q$-Gevrey regular}.
        Then $M=(M_k)$ is strongly non-quasianalytic but not of moderate growth, 
        thus not strongly regular. It is derivation closed.
        
        \subsection{\label{nmb:2.8}Spaces of $C^M$-functions}
        Let $U \subseteq \mathbb{R}^n$ be open and let $M$ be a DC-weight sequence. 
        For any $\rho >0$ and $K \subseteq U$ compact with smooth boundary, define
        \[
                C^M_\rho(K) := \{f \in C^\infty(K) : \|f\|_{\rho,K} < \infty\}
        \]
        with
        \[
                \|f\|_{\rho,K} :=
                \sup \Big\{ \frac{|\partial^{\alpha} f(x)|}{\rho^{|\alpha|} \, |\alpha|! \, M_{|\alpha|}} : \alpha \in
                \mathbb{N}^n, x \in K\Big\}.
        \]
        It is easy to see that $C^M_\rho(K)$ is a Banach space.
        In the description of $C^M_\rho(K)$, instead of compact $K$ with smooth
        boundary, we may also use open $K\subset U$ with $\overline K$ compact in
        $U$, like \cite{Thilliez08}. Or we may work with Whitney jets on compact
        $K$, like \cite{Komatsu73}.
        
        The space $C^M(U)$ carries the projective limit topology over compact $K
        \subseteq U$ of the inductive limit over $\rho \in \mathbb{N}_{>0}$:
        \[
                C^M(U) = \varprojlim_{K \subseteq U} \big(\varinjlim_{\rho \in \mathbb{N}_{>0}}
                C^M_\rho(K)\big).
        \]
        One can prove that, for $\rho < \rho'$, the canonical injection $C^M_\rho(K) \to C^M_{\rho'}(K)$
        is a compact mapping; it is even nuclear 
        (see \cite{Komatsu73}, \cite[p. 166]{Komatsu73p}).
        Hence $\varinjlim_\rho C^M_\rho(K)$ is a Silva space, i.e., an inductive limit of
        Banach spaces such that the canonical
        mappings are compact; therefore it is complete, webbed, and
        ultrabornological, see \cite{Floret71}, \cite[5.3.3]{Jarchow81}, also
        \cite[52.37]{KM97}. 
        We shall use this locally convex topology below only
        for $n=1$ -- in general it is stronger than the one which we will define in
        \thetag{\ref{nmb:3.1}}, 
        but it has the same system of bounded sets, see \thetag{\ref{nmb:4.6}}.

        \begin{lemma}\label{nmb:2.9}
                For a logarithmically convex sequence $M_k$ with $M_0=1$ we have 
                \begin{equation*}
                        M_1^k \, M_k\ge M_j\, M_{\alpha_1} \cdots M_{\alpha_j} \quad \text{ for all }
                        \alpha_i\in \mathbb{N}_{>0}, \alpha_1+\dots+\alpha_j = k.
                \end{equation*}
        \end{lemma}
        
        \demo{Proof}
        We use induction on $k$.
        The assertion is trivial for $k=j$. 
        Assume that $j < k$. Then there exists $i$ such that $\alpha_i \ge 2$.
        Put $\alpha_i':=\alpha_i-1$. By induction hypothesis,
        \[
                M_j \, M_{\alpha_1} \cdots M_{\alpha_i'} \cdots M_{\alpha_j} \le M_1^{k-1} M_{k-1}.
        \] 
        Since $M_{k+1}/M_k$ is increasing by \eqref{logconvex}, we obtain
        \begin{align*}
                M_j \, M_{\alpha_1} \cdots M_{\alpha_j} 
                &= M_j \, M_{\alpha_1} \cdots M_{\alpha_i'} \cdots M_{\alpha_j} \cdot \frac{M_{\alpha_i}}{M_{\alpha_i'}}
                \le\\&
                \le M_1^{k-1} \, M_{k-1} \cdot \frac{M_k}{M_{k-1}} \le M_1^k \, M_k.
                \qed\end{align*}
        \enddemo

        \setlength{\LTcapwidth}{4.5in}
        \begin{longtable}{|l|c|l|}
                \caption{
                        Let $M=(M_k)$ and $N=(N_k)$ be increasing $(\le)$ sequences of real numbers with $M_0=N_0=1$. 
                        By $U$ we denote an open subset of $\mathbb{R}^n$. The mapping
                        $T_a : C^M(U) \to \mathcal{F}_n^M$ is the Taylor series homomorphism for $a \in U$ (see \thetag{\ref{nmb:2.2}}).
                        Recall that $M$ is a DC-weight sequence if it is 
                        logarithmically convex and stable under derivation.} \\
                \endfirsthead
                \hline\endhead
                \hline 
                & & \\ [-.7ex]
                {\bf Properties of $M$} & & {\bf Properties of $C^M$}\\ [0.5ex]
                \hline
                \hline
                & & \\ [-1.5ex]
                $M$ increasing, $M_0=1$, & $\Rightarrow$ & $C^\omega(U) \subseteq C^M(U) \subseteq C^\infty(U)$ \\ 
                (always assumed below this line) & & \\
                [0.5ex]\hline
                & & \\ [-1.5ex]
                $M$ is logarithmically convex & $\Rightarrow$ & $C^M(U)$ is a ring.\\
                (always assumed below this line), & & $C^M$ is closed under composition.\\
                i.e., $M_k^2 \le M_{k-1}\, M_{k+1}$ for all $k$. & & $C^M$ is closed under
                applying the \\
                Then: $(M_k)^{1/k}$ is increasing, & & \qquad inverse function theorem.\\
                \hspace{.2cm} $M_l\, M_k \le M_{l+k}$ for all $l,k$, & & $C^M$ is closed under solving ODEs.\\
                \hspace{.2cm} and $M_1^k\, M_k \ge M_j\, M_{\alpha_1} \cdots M_{\alpha_j}$ & & \\
                \hspace{.2cm} for $\alpha_i \in \mathbb{N}_{>0},\alpha_1+\cdots+\alpha_j=k$. & & \\ [0.5ex]
                \hline
                & & \\ [-1.5ex]
                $\sup_{k \in \mathbb{N}_{>0}} (M_k/N_k)^{1/k} <\infty$ & $\Leftrightarrow$ & $C^M(U) \subseteq C^N(U)$\\ [0.5ex]
                \hline
                & & \\ [-1.5ex]
                $\sup_{k \in \mathbb{N}_{>0}} (M_k)^{1/k} <\infty$ & $\Leftrightarrow$ & $C^\omega(U) = C^M(U)$\\ [0.5ex]
                \hline
                & & \\ [-1.5ex]
                $\lim_{k \to \infty} (M_k)^{1/k} =\infty$ & $\Leftrightarrow$ & $C^\omega(U) \subsetneq C^M(U)$\\ [0.5ex]
                \hline
                & & \\ [-1.5ex]
                $\sup_{k \in \mathbb{N}_{>0}} (M_{k+1}/M_k)^{1/k} <\infty$ & $\Leftrightarrow$ & $C^M$ is closed under derivation.\\
                (always assumed below this line) & & \\ [0.5ex]
                \hline
                & & \\ [-1.5ex]
                $\sum_{k=0}^\infty \frac{M_k}{(k+1)M_{k+1}} =\infty$ & $\Leftrightarrow$ & $C^M$ is quasianalytic,\\
                \qquad or, equivalently, & & i.e., $T_a : C^M(U) \to \mathcal{F}_n^M$ is injective\\ 
                $\sum_{k=1}^\infty (\frac1{k!\,M_k})^{1/k} = \infty$ & & 
                (not surjective if $C^\omega(U) \subsetneq C^M(U)$).\\ [0.5ex]
                \hline
                & & \\ [-1.5ex]
                $\sum_{k=0}^\infty \frac{M_k}{(k+1)M_{k+1}} <\infty$ & $\Leftrightarrow$ & $C^M$ is non-quasianalytic.\\ 
                & & Then $C^M$ partitions of unity exist.\\ [0.5ex]
                \hline
                & & \\ [-1.5ex]
                $\lim_{k \to \infty} (M_k)^{1/k} =\infty$ and & $\Leftrightarrow$ & $C^\omega(U) \subsetneq C^M(U)$ and\\ 
                $\sum_{k=j}^\infty \frac{M_k}{(k+1)M_{k+1}} \le C \frac{M_j}{M_{j+1}}$ & & 
                $T_a : C^M(U) \to \mathcal{F}_n^M$ is surjective, i.e.,\\ 
                for all $j \in\mathbb{N}$ and some $C$ & & $C^M$ is strongly non-quasianalytic.\\[0.5ex]
                \hline
                & & \\ [-1.5ex]
                $M$ has moderate growth, i.e., & $\Rightarrow$ & $C^M$ is cartesian closed \\
                $\sup_{j,k \in \mathbb{N}_{>0}} (\frac{M_{j+k}}{M_j\, M_k})^{1/(j+k)} < \infty$ & & 
                will be proved in \thetag{\ref{nmb:5.3}} \\ [0.5ex]
                \hline
                & & \\ [-1.5ex]
                $M$ is strongly regular, i.e., & $\Rightarrow$ & Whitney's extension theorem\\
                it is strongly non-quasianalytic & & holds in $C^M$.\\ 
                and has moderate growth. & & \\[0.5ex]
                \hline
                & & \\ [-1.5ex]
                $\delta>0$ and $M_k = (k!)^\delta$ for $k \in \mathbb{N}$. &$\Leftrightarrow$ & $C^M$ is the Gevrey class $G^{1+\delta}$.\\
                Then $M$ is strongly regular. & & \\ [0.5ex]
                \hline 
        \end{longtable}
        
        \section{\label{nmb:3} $C^M$-mappings}
        
        \subsection{\label{nmb:3.1}Definition: $C^M$-mappings}
        Let $M$ be a DC-weight sequence, and 
        let $E$ be a locally convex vector space. 
        A curve $c:\mathbb R\to E$ is called $C^M$ if for each continuous linear
        functional $\ell\in E^*$ the curve $\ell\circ c:\mathbb R\to \mathbb R$ is of
        class $C^M$. 
        The curve $c$ is called {\it strongly $C^M$} if $c$ is smooth and 
        for all
        compact $K \subset \mathbb R$ 
        there exists $\rho>0$ such that
        \[
        \left\{\frac{c^{(k)}(x)}{ \rho^{k} \, k! \, M_{k}}: 
        k \in \mathbb{N}, x \in K\right\}\text{ is bounded in }E.
        \]
        The curve $c$ is called {\it strongly uniformly $C^M$} if $c$ is smooth and 
        there exists $\rho>0$ such that
        \[
        \left\{\frac{c^{(k)}(x)}{ \rho^{k} \, k! \, M_{k}}: 
        k \in \mathbb{N}, x \in \mathbb R\right\}\text{ is bounded in }E.
        \]
        
        Now let $M$ be a non-quasianalytic DC-weight sequence.
        Let $U$ be a $c^\infty$-open subset of $E$, and let $F$ be another 
        locally convex vector space. 
        {\it A mapping $f:U\to F$ is called $C^M$ if $f$ is smooth 
                in the sense of \thetag{\ref{nmb:7.3}}
                and if $f\circ c$ is a 
                $C^M$-curve in $F$ for
                every $C^M$-curve $c$ in $U$}.
        Obviously, 
        {\it the composite of $C^M$-mappings is again a $C^M$-mapping, and the chain rule holds}.
        This notion is equivalent to the expected one on Banach spaces, see \ref{nmb:3.9} below.
        
        We equip the space $C^M(U,F)$ with the initial locally convex structure
        with respect to the family of mappings 
        \begin{gather*}
                C^M(U,F) \East{C^M(c,\ell)}{} C^M(\mathbb R,\mathbb R),
                \quad f\mapsto \ell\circ f\circ c,\quad \ell \in E^*, c\in C^M(\mathbb R,U)
        \end{gather*}
        where $C^M(\mathbb R,\mathbb R)$ carries the locally convex structure
        described in \thetag{\ref{nmb:2.8}} and where $E^*$ is the space of all continuous linear functionals 
        on $E$.
        
        For $U\subseteq \mathbb R^n$, 
        this locally convex topology differs from the one described in \thetag{\ref{nmb:2.8}},
        but they have the same bounded sets, see \thetag{\ref{nmb:4.6}} below. 
        
        If $F$ is convenient, then
        by standard arguments, %
        the space $C^M(U,F)$ is $c^\infty$-closed in the product 
        $\prod_{\ell,c}C^M(\mathbb R,\mathbb R)$ 
        and hence is {\it convenient}. 
        If $F$ is convenient, then a mapping $f:U\to F$ is $C^M$ if and only if $\ell\circ f$ is $C^M$ for all 
        $\ell\in F^*$.

        \subsection{\label{nmb:3.2} Example: There are weak $C^M$-curves which are not strong}
        By \cite[Theorem 1]{Thilliez08}, for each DC-weight sequence $M$ there exists 
        $f \in C^M(\mathbb{R},\mathbb{R})$ such that $|f^{(k)}(0)| \ge k! \, M_k$ for all $k \in \mathbb{N}$. 
        Then $g:\mathbb R\to \mathbb R^{\mathbb{N}}$ given by 
        $g(t)_n = f(nt)$ is $C^M$ but not strongly $C^M$.
        Namely, each bounded linear functional $\ell$ on $\mathbb R^{\mathbb{N}}$
        depends only on finitely many coordinates, so we take the maximal $\rho$ for the finitely many 
        coordinates of $g$ being involved.
        On the other hand, for each $\rho$ and any compact neighborhood $L$ of $0$ 
        the set 
        \[
        \left\{\frac{g^{(k)}(t)}{\rho^{k}\, k!\, M_{k}}:t\in L, k\in \mathbb{N}\right\}
        \]
        has $n$-th coordinate unbounded if $n>\rho$.
        
        \begin{lemma}\label{nmb:3.3}
                Let $E$ be a convenient vector space such that there exists a Baire vector
                space topology on the dual $E^*$ for which the point evaluations $\on{ev}_x$
                are continuous for all $x\in E$.
                Then a curve $c:\mathbb R\to E$ is $C^M$ if and
                only if $c$ is strongly $C^M$, for any DC-weight sequence $M$. 
        \end{lemma}
        
        See \thetag{\ref{nmb:5.2}} for a more general version.
        
        \begin{demo}{Proof}
                Let $K$ be compact in $\mathbb R$. 
                We consider the sets
                \[
                A_{\rho,C} :=\Bigl\{\ell\in E^*: \frac{|(\ell\circ c)^{(k)}(x)|}{\rho^{k}\,
                        k!\, M_{k}}\le C\text{ for all }k\in \mathbb{N}, x\in K\Bigr\}
                \]
                which are closed subsets in $E^*$ for the Baire topology. We have
                $\bigcup_{\rho,C}A_{\rho,C}= E^*$. By the Baire property there exists $\rho$
                and $C$ such that the interior $U$ of $A_{\rho,C}$ is non-empty. If
                $\ell_0\in U$ then for all $\ell\in E^*$ there is an $\epsilon>0$ such that 
                $\epsilon\ell\in U-\ell_0$ and hence for all $x\in K$ and all $k$ we have
                \begin{align*}
                        |(\ell\circ c)^{(k)}(x)| \le \tfrac1\epsilon \left(|((\epsilon\ell+\ell_0)\circ c)^{(k)}(x)| +
                        |(\ell_0\circ c)^{(k)}(x)|\right) \le \tfrac{2C}{\epsilon}\,\rho^{k}\,k!\,
                        M_{k}.
                \end{align*}
                So the set 
                \[
                \left\{\frac{c^{(k)}(x)}{\rho^{k}\, k!\, M_{k}}: k\in\mathbb{N}, x\in K
                \right\}
                \]
                is weakly bounded in $E$ and hence bounded. 
                \qed\end{demo}

        \begin{lemma}\label{nmb:3.4}%
                Let $M$ be a DC-weight sequence, and let $E$ be a Banach space.
                For a curve $c:\mathbb R\to E$ the following are equivalent.
                \begin{enumerate}
                \item $c$ is $C^M$.
                \item For each sequence $(r_k)$ with $r_k\,t^k\to
                        0$ for all $t>0$, and each compact set $K$ in $\mathbb R$, 
                        the set $\{\frac1{k!M_k}\,c^{(k)}(a)\,r_k: a\in K, k\in \mathbb{N}\}$ is 
                        bounded in $E$.
                \item For each sequence $(r_k)$ satisfying $r_k>0$, 
                        $r_kr_\ell\geq r_{k+\ell}$, 
                        and $r_k\,t^k\to 0$ for all $t>0$, and each compact set $K$ in $\mathbb R$, 
                        there exists an
                        $\epsilon>0$ such that 
                        $\{\frac1{k!M_k}\,c^{(k)}(a)\,r_k\,\epsilon^k: a\in K, k\in \mathbb{N}\}$ 
                        is bounded in $E$. 
                \end{enumerate}
        \end{lemma}
        
        \demo{Proof}
        {\rm (1)} $\implies$ {\rm (2)}
        For $K$, there exists $\rho>0$ such that
        \[
        \left\|\frac{c^{(k)}(a)}{k!\,M_k}r_k \right\|_E = 
        \left\|\frac{c^{(k)}(a)}{k!\,\rho^k\,M_k}\right\|_E\cdot |r_k\rho^k| 
        \]
        is bounded uniformly in $k\in \mathbb N$ and $a\in K$ by \thetag{\ref{nmb:3.3}}. 
        
        {\rm (2)} $\implies$ {\rm (3)} Use $\epsilon=1$.
        
        {\rm (3)} $\implies$ {\rm (1)} 
        Let $a_k:=\sup_{a\in K}\|\frac1{k!\,M_k}\,c^{(k)}(a)\|_E$. Using 
        \cite[9.2.(4$\Rightarrow$1)]{KM97} 
        these are the coefficients of a power series with positive 
        radius of convergence.
        Thus $a_k/\rho^k$ is bounded for some $\rho>0$. 
        \qed\enddemo

        \begin{lemma}\label{nmb:3.5}
                Let $M$ be a DC-weight sequence. 
                Let $E$ be a convenient vector space, and let $\mathcal S$ be a
                family of bounded linear functionals on $E$ which together detect bounded
                sets (i.e., $B\subseteq E$ is bounded if and only if $\ell(B)$ is bounded for all
                $\ell\in\mathcal S$). 
                Then a curve $c:\mathbb R\to E$ is $C^M$ if and only if 
                $\ell\circ c:\mathbb R\to \mathbb R$ is $C^M$
                for all $\ell\in\mathcal S$.
        \end{lemma}
        
        \begin{demo}{Proof}
                For smooth curves this follows from \cite[2.1 and 2.11]{KM97}.
                By \thetag{\ref{nmb:3.4}}, for any $\ell\in E'$, 
                the function $\ell\circ c$ is $C^M$ if and only if: 
                \begin{enumerate}
                \item
                        For each sequence $(r_k)$ with $r_k\,t^k\to
                        0$ for all $t>0$, and each compact set $K$ in $\mathbb R$, 
                        the set $\{\frac1{k!M_k}\,(\ell\circ c)^{(k)}(a)\,r_k: a\in K, k\in \mathbb{N}\}$ is 
                        bounded. 
                \end{enumerate}
                By (1) the curve $c$ is $C^M$ if and only if 
                the set $\{\frac1{k!M_k}\,c^{(k)}(a)\,r_k: a\in K, k\in \mathbb{N}\}$ is
                bounded in $E$. By (1) again this is in turn equivalent to 
                $\ell\circ c\in C^M$ for all $\ell\in\mathcal S$, 
                since $\mathcal S$ detects bounded sets.
                \qed\end{demo}

        \subsection{\label{nmb:3.6} $C^M$ curve lemma}
        A sequence $x_n$ in a locally convex space $E$ 
        is said to be
        {\it Mackey convergent} to $x$, if there exists some $\lambda_n\nearrow \infty$
        such that $\lambda_n(x_n-x)$ is bounded. If we fix $\lambda=(\lambda_n)$ we say that 
        $x_n$ is $\lambda$-converging.

        \begin{lemma*} 
                Let $M$ be a non-quasianalytic DC-weight sequence.
                Then there exist sequences $\lambda_k\to 0$, $t_k\to t_\infty$, $s_k>0$ in $\mathbb R$ 
                with the following property:
                For $1/\lambda=(1/\lambda_n)$-converging sequences $x_n$ and $v_n$
                in a convenient vector space $E$ 
                there exists a strongly uniformly $C^M$-curve $c:\mathbb R\to E$ with
                $c(t_k+t)= x_k+t.v_k$ for $|t|\leq s_k$.
        \end{lemma*}
        
        \begin{demo}{Proof}
                Since $C^M$ is not quasianalytic we have $\sum_k 1/(k!M_k)^{1/k}<\infty$. We choose
                another non-quasianalytic DC-weight sequence $\bar M= (\bar M_k)$ with 
                $(M_k/\bar M_k)^{1/k}\to \infty$. By \thetag{\ref{nmb:2.3}} there is a $C^{\bar M}$-function 
                $\varphi:\mathbb R\to [0,1]$ which is 0 on $\{t: |t|\ge \frac12\}$ and which is 1 on
                $\{t: |t|\le \frac13\}$, i.e.\  there exist $\bar C,\rho>0$ such that
                \[
                        |\varphi^{(k)}(t)|\leq \bar C\,\rho^k\,k!\,\bar M_k
                        \quad\text{ for all }t\in\mathbb R\text{ and }k\in\mathbb{N}.
                \]
                For $x,v$ in a absolutely convex bounded set $B\subseteq E$ and $0<T\le 1$ the curve
                $c:t\mapsto \varphi(t/T)\cdot (x+t\,v)$ satisfies
                (cf.\ \cite[Lemma 2]{Boman67}):
                \begin{align*}
                        c^{(k)}(t)&= T^{-k}\varphi^{(k)}(\tfrac{t}{T}).(x+t.v)+ k\,T^{1-k}\,\varphi^{(k-1)}(\tfrac{t}{T}).v
                        \\&
                        \in T^{-k}\bar C\,\rho^k\,k!\,\bar M_k(1+\tfrac{T}{2}).B+
                        k\,T^{1-k}\,\bar C\,\rho^{k-1}\,(k-1)!\,\bar M_{k-1}.B
                        \\&
                        \subseteq T^{-k}\bar C\,\rho^k\,k!\,\bar M_k(1+\tfrac{T}{2}).B+
                        T\,T^{-k}\,\bar C\tfrac1\rho\,\rho^{k}\,k!\,\bar M_{k}.B
                        \\&
                        \subseteq \bar C(\tfrac32+ \tfrac 1\rho)\,T^{-k}\,\rho^{k}\,k!\,\bar M_{k}.B
                \end{align*}
                So there are $\rho,C:=\bar C(\tfrac32+ \tfrac 1\rho)>0$ which do not depend on $x,v$ and $T$ such that
                $c^{(k)}(t)\in C\, T^{-k}\,\rho^k\,k!\,\bar M_k. B$ for all $k$ and $t$.
                
                Let $0<T_j\leq 1$ with $\sum_j T_j<\infty$
                and $t_k:=2\sum_{j<k}T_j+T_k$.
                We choose the $\lambda_j$ such 
                that
                $0<\lambda_j/T_j^k\leq M_k/\bar M_k$ (note that $T_j^k\,M_k/\bar M_k\to\infty$ for $k\to\infty$) 
                for all $j$ and $k$, and that
                $\lambda_j/T_j^k\to 0$ for $j\to\infty$ and each $k$.
                
                Without loss we may assume that $x_n\to 0$. 
                By assumption there exists a closed bounded absolutely convex subset $B$ in $E$ 
                such that 
                $x_n,v_n\in \lambda_n\cdot B$. 
                We consider
                $c_j:t\mapsto \varphi\bigl((t-t_j)/T_j\bigr)\cdot \bigl(x_j+(t-t_j)\,v_j\bigr)$
                and $c:=\sum_j c_j$.
                The $c_j$ have disjoint support $\subseteq [t_j-T_j,t_j+T_j]$, hence $c$ is $C^\infty$ on 
                $\mathbb R\setminus\{t_\infty\}$ with
                \[
                        c^{(k)}(t) \in C\, T_j^{-k} \,\rho^{k} k!\bar M_k\,\lambda_j\cdot B 
                        \quad\text{ for }|t-t_j|\leq T_j.
                \]
                Then
                \[
                        \|c^{(k)}(t)\|_{B}\le C\,\rho^{k}\,k!\bar M_k\,\frac{\lambda_j}{T_j^{k}} 
                        \leq C \rho^{k} k!\bar M_k\,\frac{M_k}{\bar M_k} = C\,\rho^{k}\,k! M_k
                \]
                for $t\ne t_\infty$.
                Hence $c:\mathbb R\to E_B$ 
                (see \cite[2.14.6]{KM97} or \thetag{\ref{nmb:7.1}})
                is smooth at $t_\infty$ as well, 
                and is strongly $C^M$ by
                the following lemma.%
                \qed\end{demo}
        
        \begin{lemma}\label{nmb:3.7}
                Let $c:\mathbb R\setminus\{0\}\to E$ be strongly $C^M$ in the sense that 
                $c$ is smooth and for all
                bounded $K \subset \mathbb R\setminus \{0\}$ 
                there exists $\rho>0$ such that
                \[
                \left\{\frac{c^{(k)}(x)}{ \rho^{k} \, k! \, M_{k}}: 
                k \in \mathbb{N}, x \in K\right\}\text{ is bounded in }E.
                \]
                Then $c$ has a unique extension to a strongly $C^M$-curve on $\mathbb R$.
        \end{lemma}
        
        \begin{demo}{Proof}
                The curve $c$ has a unique extension to a smooth curve by \cite[2.9]{KM97}.
                The strong $C^M$ condition extends by continuity. 
                \qed\end{demo}
        
        \begin{corollary}\label{nmb:3.8}
                Let $M$ be a non-quasianalytic DC-weight sequence. Then we have:
                \begin{enumerate}
                \item
                        The final topology on $E$ with respect to all strongly $C^M$-curves equals the
                        Mackey closure topology.
                \item
                        A locally convex space $E$ is convenient \thetag{\ref{nmb:7.2}} if and only if 
                        for any (strongly) $C^M$-curve $c:\mathbb R\to E$ there exists a (strongly) 
                        $C^M$-curve
                        $c_1:\mathbb R\to E$ with $c_1'=c$.
                \end{enumerate}
        \end{corollary}
        \begin{demo}{Proof}
                (1) For any Mackey converging sequence there exists a $C^M$-curve passing 
                through a subsequence in finite time by \thetag{\ref{nmb:3.6}}. So the final topologies
                generated by the Mackey converging sequences and by the $C^M$-curves coincide. 
                
                (2) In order to show that a locally convex space $E$ is convenient, we 
                have to prove that it is $c^\infty$-closed in its completion. So let
                $x_n\in E$ converge Mackey to $x_\infty$ in the completion.
                Then by \thetag{\ref{nmb:3.6}} there exists a strongly $C^M$-curve $c$ in the completion
                passing in finite time through a subsequence of the $x_n$ with velocity $v_n=0$. 
                The form of $c$
                (in the proof of \thetag{\ref{nmb:3.6}})
                shows that its derivatives $c^{(k)}(t)$ for $k>0$ are multiples of the $x_n$ and hence 
                have values in $E$. Then $c'$ is a $C^M$-curve and so the antiderivative $c$ of $c'$
                lies in $E$ by assumption. In particular $x_\infty\in c(\mathbb R)\subseteq E$.
                
                Conversely, if $E$ is convenient, then every smooth curve $c$ has a smooth 
                antiderivative $c_1$ in $E$ by \cite[2.14]{KM97}. Since
                \[
                        \frac1{\rho^{k+1}(k+1)!\,M_{k+1}}c_1^{(k+1)}(t)=
                        \frac{M_k}{\rho(k+1) M_{k+1}}\frac1{\rho^k\,k!\,M_k}c^{(k)}(t)
                \]
                and since
                \[
                        \frac{M_k}{\rho(k+1)M_{k+1}}\leq \frac1{\rho M_1} 
                \]
                by \thetag{\ref{nmb:2.1}.2} the antiderivative $c_1$ is (strongly) $C^M$ if $c$ is so.
                \qed\end{demo}
        
        \begin{theorem}\label{nmb:3.9}
                Let $M=(M_k)$ be a non-quasianalytic DC-weight sequence.
                Let $U\subseteq E$ be $c^\infty$-open in a convenient vector space, and let $F$
                be a Banach space. 
                For a mapping $f:U\to F$, the following assertions are
                equivalent.
                \begin{enumerate}
                \item $f$ is $C^M$. 
                \item $f$ is $C^M$ along strongly $C^M$ curves.
                \item $f$ is smooth, and for each closed bounded absolutely convex $B$ in $E$
                        and each $x\in U\cap E_B$ there are $r>0$, $\rho>0$, and $C>0$ such that
                        \[
                        \frac1{k!\,M_k}\|d^k(f\circ i_B)(a)\|_{L^k(E_B,F)} \le C\,\rho^k
                        \]
                        for all $a\in U\cap E_B$ with $\|a-x\|_{B}\le r$ and all $k\in\mathbb{N}$.
                \item $f$ is smooth, and for each closed bounded absolutely convex $B$ in $E$
                        and each compact $K\subseteq U\cap E_B$ there are $\rho>0$ and $C>0$ such that
                        \[
                        \frac1{k!\,M_k}\|d^k(f\circ i_B)(a)\|_{L^k(E_B,F)} \le C\,\rho^k
                        \]
                        for all $a\in K$ and all $k\in\mathbb{N}$.
                \end{enumerate}
        \end{theorem}
        
        \begin{demo}{Proof}
                {\rm (1)} $\implies$ {\rm (2)} is clear.
                
                {\rm (2)} $\implies$ {\rm (3)}
                Without loss let $E=E_B$ be a Banach space. 
                For each $v\in E$ and $x\in U$ the iterated directional derivative 
                $d_v^kf(x)$ exists since $f$ is $C^M$ along affine lines. 
                To show that $f$ is smooth it suffices to check that
                $d_{v_n}^kf(x_n)$ is bounded for each $k\in \mathbb{N}$ and each Mackey
                convergent sequences $x_n$ and $v_n\to 0$, by \cite[5.20]{KM97}. 
                For contradiction let us assume that there exist $k$ and sequences $x_n$
                and $v_n$ with $\|d_{v_n}^kf(x_n)\|\to \infty$. 
                By passing to a subsequence 
                we may assume that $x_n$ and $v_n$ are $(1/\lambda_n)$-converging
                for the $\lambda_n$ from \thetag{\ref{nmb:3.6}}.
                Hence there exists a strongly $C^M$-curve $c$
                in $E$ and with $c(t+t_n)=x_{n}+t.v_{n}$ for $t$ near 0 for each
                $n$ separately, and for $t_n$ from \thetag{\ref{nmb:3.6}}. 
                But then 
                $\|(f\circ c)^{(k)}(t_n)\|=\|d_{v_{n}}^k f(x_{n})\|\to \infty$, a contradiction.
                So $f$ is smooth. 
                
                Assume for contradiction that the boundedness condition in
                {\rm (3)} does not hold. Then there exists $x\in U$ such that for
                all $r,\rho, C>0$ there is an $a=a(r,\rho,C)\in U$ and $k=k(r,\rho,C)\in \mathbb{N}$ with
                $\|a-x\|\le r$ but 
                \[
                \frac1{k!\,M_k}\|d^kf(a)\|_{L^k(E,F)} > C\,\rho^k.
                \]
                By 
                \cite[7.13]{KM97}
                we have
                \[
                \|d^kf(a)\|_{L^k(E,F)} %
                \le (2e)^k\sup_{\|v\|\le 1}\|d_v^kf(a)\|.
                \]
                So for each $\rho$ and $n$ take $r=\tfrac1{n\rho}$ and $C=n$. Then there are
                $a_{n,\rho}\in U$ with $\|a_{n,\rho}- x\|\le \tfrac1{n\rho}$, moreover
                $v_{n,\rho}$ with $\|v_{n,\rho}\|=1$,
                and $k_{n,\rho}\in \mathbb{N}$ such that
                \[
                \frac{(2e)^{k_{n,\rho}}}{k_{n,\rho}!\,M_{k_{n,\rho}}\,\rho^{k_{n,\rho}}}
                \|d_{v_{n,\rho}}^{k_{n,\rho}}f(a_{n,\rho})\| > n.
                \]
                Since $K:=\{a_{n,\rho}:n,\rho\in \mathbb{N} \}\cup\{x\}$ is compact, this
                contradicts the following
                
                {\bf Claim.} {\it For each compact $K\subseteq E$ there are $C,\rho\geq 0$ such that for all
                        $k\in\mathbb{N}$ and $x\in K$ we have 
                        $\sup_{\|v\|\le 1}\|d_v^k f (x)\|\leq C\,\rho^{k} k! M_k$.}
                \newline
                Otherwise, there exists a compact set $K\subseteq E$ such that 
                for each $n\in\mathbb{N}$ there are $k_n\in\mathbb{N}$, $x_n\in K$, and $v_n$ with
                $\|v_n\|=1$ such that
                \[
                        \|d_{v_n}^{k_n}f(x_n)\|>{k_n}!\,M_{k_n}\,\left(\frac1{\lambda_n^2}\right)^{k_n+1},
                \]
                where we used $C=\rho:=1/\lambda_n^2$ with the $\lambda_n$ from \thetag{\ref{nmb:3.6}}.
                By passing to a subsequence (again denoted $n$) we may assume that the $x_n$ are
                $1/\lambda$-converging, thus
                there exists a strongly $C^M$-curve $c:\mathbb R\to E$ 
                with $c(t_n+t)=x_n+t.\lambda_n.v_n$ for $t$ near 0 by \thetag{\ref{nmb:3.6}}.
                Since
                \[
                        (f\circ c)^{(k)}(t_n)= \lambda_n^k d_{v_n}^k f(x_n),
                \]
                we get 
                \[
                        \left(\frac{\|(f\circ c)^{(k_n)}(t_n)\|}{{k_n}!M_{k_n}}\right)^{\frac1{k_n+1}}
                        =\left(\lambda_n^{k_n} \frac{\|d_{v_n}^{k_n}f(x_n)\|}{{k_n}!M_{k_n}}\right)^{\frac1{k_n+1}}
                        > \frac1{\lambda_n^{\frac{k_n+2}{k_n+1}}} \to\infty,
                \]
                a contradiction to $f\circ c\in C^M$.
                
                {\rm (3)} $\implies$ {\rm (4)} is obvious since the compact
                set $K$ is covered by finitely many balls. 
                
                {\rm (4)} $\implies$ {\rm (1)}
                We have to show that $f\circ c$ is $C^M$ for each $C^M$-curve $c:\mathbb R\to
                E$. By \thetag{\ref{nmb:3.4}.2} it suffices to show that for each 
                sequence $(r_k)$ satisfying $r_k>0$, 
                $r_kr_\ell\geq r_{k+\ell}$, 
                and $r_k\,t^k\to 0$ for all $t>0$, and each compact interval $I$ in $\mathbb R$, 
                there exists an
                $\epsilon>0$ such that 
                $\{\frac1{k!M_k}\,(f\circ c)^{(k)}(a)\,r_k\,\epsilon^k: a\in I, k\in \mathbb{N}\}$ 
                is bounded.
                
                By \thetag{\ref{nmb:3.4}.2} applied to $r_k2^k$ instead of $r_k$, for each $\ell\in E^*$, 
                each sequence $(r_k)$ with $r_k\,t^k\to
                0$ for all $t>0$, and each compact interval $I$ in $\mathbb R$
                the set 
                $\{\frac1{k!M_k}\,(\ell\circ c)^{(k)}(a)\,r_k\,2^k: a\in I, k\in \mathbb{N}\}$ is 
                bounded in $\mathbb R$.
                Thus $\{\frac1{k!M_k}\,c^{(k)}(a)\,r_k\,2^k: a\in I, k\in \mathbb{N}\}$ 
                is contained in some closed absolutely convex $B\subseteq E$. 
                Consequently, $c^{(k)}:I\to E_B$ is smooth and hence
                $K_k:=\{\frac1{k!M_k}\,c^{(k)}(a)\,r_k\,2^k: a\in I\}$ is
                compact in $E_B$ for each $k$. 
                Then each sequence $(x_n)$ in the set
                \[
                K:=\left\{\frac1{k!M_k}\,c^{(k)}(a)\,r_k: a\in I, k\in \mathbb{N}\right\}
                =\bigcup_{k\in\mathbb{N}}\frac1{2^k} K_k
                \]
                has a cluster point in $K\cup\{0\}$: either there is a subsequence in one
                $K_k$, or $2^{k_n}x_{k_n}\in K_{k_n}\subseteq B$ for $k_n\to\infty$, hence
                $x_{k_n}\to 0$ in $E_B$. 
                So $K\cup\{0\}$ is compact. 
                
                By Fa\`a di Bruno (\cite{FaadiBruno1855} for the 1-dimensional version) 
                \begin{align*}
                        \frac{(f\circ c)^{(k)}(a)}{k!} = 
                        \sum_{j\ge 0} \sum_{\substack{\alpha\in \mathbb{N}_{>0}^j\\ \alpha_1+\dots+\alpha_j =k}}
                        \frac{1}{j!}d^jf(c(a))\Big( 
                        \frac{c^{(\alpha_1)}(a)}{\alpha_1!},\dots,
                        \frac{c^{(\alpha_j)}(a)}{\alpha_j!}\Big)
                \end{align*}
                and \eqref{logconvex2} for $a\in I$ and $k\in \mathbb{N}$ we have
                \begin{align*}
                        &\left\|\frac1{k!M_k}\,(f\circ c)^{(k)}(a)\,r_k \right\|\le
                        \\&
                        \le M_1^k \sum_{j\ge 0} \sum_{\substack{\alpha\in \mathbb{N}_{>0}^j\\ \alpha_1+\dots+\alpha_j =k}}
                        \frac{\|d^jf(c(a))\|_{L^j(E_B,F)}}{j!M_j}\prod_{i=1}^j
                        \frac{\|c^{(\alpha_i)}(a)\|_B \, r_{\alpha_i}}{\alpha_i!M_{\alpha_i}}
                        \\&
                        \le M_1^k \sum_{j\ge 0} \binom{k-1}{j-1} C\,\rho^j\, \frac1{2^k}
                        = M_1^k \rho(1+\rho)^{k-1} C\, \frac1{2^k}.
                \end{align*}
                So 
                $
                \left\{\frac1{k! M_k}\,(f\circ c)^{(k)}(a)\,\Big(\frac{2}{M_1(1+\rho)}\Big)^k\,r_k: 
                a\in I, k\in \mathbb{N}\right\}
                $
                is bounded as required.
                \qed\end{demo}
        
        \begin{corollary}\label{nmb:3.10}
                Let $M$ and $N$ be non-quasianalytic DC-weight sequences with \eqref{incl}
                \begin{equation*} 
                        \sup_{k \in \mathbb{N}_{>0}} \Big(\frac{M_k}{N_k}\Big)^{\frac{1}{k}} < \infty.
                \end{equation*}
                Then $C^M(U,F)\subseteq C^N(U,F)$ for all convenient vector spaces $E$ and
                $F$ and each $c^\infty$-open $U\subseteq E$. 
                Moreover $C^\omega(U,F)\subseteq C^M(U,F)\subseteq C^\infty(U,F)$.
                All these inclusions are bounded.
        \end{corollary}
        
        \begin{demo}{Proof} The inclusions
                $C^M\subseteq C^N\subseteq C^\infty$ follow from \thetag{\ref{nmb:3.9}} 
                since this is true for condition \thetag{\ref{nmb:3.9}.3} applied to $\ell\circ f$ for $\ell\in F^*$.
                
                Without loss let $F=\mathbb R$.
                If $f$ is $C^\omega$ then 
                for each closed absolutely convex bounded $B\subseteq E$ the mapping 
                $f\circ i_B:U\cap E_B\to \mathbb R$ is given by its locally converging Taylor
                series by \cite[10.1]{KM97}. So \thetag{\ref{nmb:3.9}.3} is satisfied for
                $M_k=1$ and thus for each DC-weight sequence $M$. So $f$ is $C^M$. 
                All inclusions are bounded by the uniform boundedness principle \ref{nmb:4.1} below for $C^M$ and 
                \cite[5.26]{KM97} for $C^\infty$.
                \qed\end{demo}
        
        \begin{corollary}\label{nmb:3.11}
                Let $M=(M_k)$ be a non-quasianalytic DC-weight sequence. Then we have:
                \begin{enumerate}
                \item Multilinear mappings between convenient vector spaces are $C^M$ if and only
                        if they are bounded.
                \item 
                        If $f:E\supseteq U\to F$ is $C^M$, then the derivative 
                        $df:U\to L(E,F)$ is $C^M$, and also $\widehat{df}:U\times E\to F$ is 
                        $C^M$, 
                        where the space $L(E,F)$ of all bounded linear mappings is considered 
                        with the topology of uniform convergence on bounded sets.
                \item The chain rule holds.
                \end{enumerate}
        \end{corollary}
        
        \begin{demo}{Proof}
                {\rm (1)}
                If $f$ is multilinear and $C^M$ then it is smooth by \thetag{\ref{nmb:3.9}} and hence
                bounded by \thetag{\ref{nmb:7.3}.2}. Conversely, if $f$ is multilinear and bounded then
                it is smooth by \thetag{\ref{nmb:7.3}.2}. 
                Furthermore, $f\circ i_B$ is multilinear and continuous and all derivatives of high order
                vanish. Thus condition \thetag{\ref{nmb:3.9}.3} is satisfied, so $f$ is $C^M$. 
                
                {\rm (2)} Since $f$ is smooth, by \thetag{\ref{nmb:7.3}.3} 
                the map $df:U\to L(E,F)$ exists and is smooth. 
                Let $c:\mathbb R\to U$ be a $C^M$-curve. We have to show 
                that $t\mapsto df(c(t))\in L(E,F)$ is $C^M$. By \cite[5.18]{KM97} and
                \thetag{\ref{nmb:3.5}} it suffices to show that 
                $t\mapsto c(t)\mapsto \ell(df(c(t)).v)\in \mathbb R$ is $C^M$ for each
                $\ell\in F^*$ and $v\in E$.
                We are reduced to show that $x\mapsto \ell(df(x).v)$ satisfies the
                conditions of \thetag{\ref{nmb:3.9}}. 
                By \thetag{\ref{nmb:3.9}} applied to $\ell\circ f$, 
                for each closed bounded absolutely convex $B$ in $E$
                and each $x\in U\cap E_B$ there are $r>0$, $\rho>0$, and $C>0$ such that
                \[
                \frac1{k!\,M_k}\|d^k(\ell\circ f\circ i_B)(a)\|_{L^k(E_B,\mathbb R)} \le C\,\rho^k
                \]
                for all $a\in U\cap E_B$ with $\|a-x\|_{B}\le r$ and all $k\in\mathbb{N}$. 
                For $v\in E$ and those $B$ containing $v$ we then have
                \begin{align*}
                        \|d^k(d(\ell\circ f)(&\quad)(v))\circ i_B)(a)\|_{L^k(E_B,\mathbb R)} 
                        =\|d^{k+1}(\ell\circ f\circ i_B)(a)(v,\dots)\|_{L^k(E_B,\mathbb R)} 
                        \\&
                        \le\|d^{k+1}(\ell\circ f\circ i_B)(a)\|_{L^{k+1}(E_B,\mathbb R)}\|v\|_{E_B} 
                        \le C\,\rho^{k+1}\, (k+1)!\, M_{k+1}
                        \\&
                        \le C\,\rho^k k! M_k \Big((k+1)\rho \frac{M_{k+1}}{M_k} \Big)
                        \\&
                        \le C\,\bar\rho^k k! M_k\quad \text{ for }
                        \bar \rho> \rho\; \sup_{k\geq 1}\Big((k+1)\rho \frac{M_{k+1}}{M_k} \Big)^{1/k},
                \end{align*}
                the latter quantity being finite by \eqref{der}.
                By \thetag{\ref{nmb:4.2}} below also $\widehat{df}$ is $C^M$.
                
                {\rm (3)} This is valid for all smooth $f$. 
                \qed\end{demo}
        
        \subsection{\label{nmb:3.12} Remark} For a quasianalytic DC-weight sequence $M$
        Theorem \ref{nmb:3.9} is wrong. In fact, take any rational function, e.g.
        $\frac{xy^2}{x^2+y^2}$. Let $t\mapsto x(t),y(t)$ be in $C^M(\mathbb R,\mathbb R)$ 
        with $x(0)=0=y(0)$. 
        Then $x(t)=t^r\bar x(t)$ and $y(t)=t^r\bar y(t)$ for $r>0$ and for $C^M$-functions 
        $\bar x$ and $\bar y$ since
        $C^M$ is derivation closed. 
        If $(x,y)$ is not constant we may choose $r$ such that 
        $\bar x(0)^2 + \bar y(0)^2\ne 0$, since 
        $C^M$ is quasianalytic. Then 
        $t\mapsto \frac{x(t)y(t)^2}{x(t)^2+y(t)^2}
        =t^r\frac{\bar x(t)\bar y(t)^2}{\bar x(t)^2+\bar y(t)^2}$
        is $C^M$ near 0, but the rational function is not smooth.

        \section{\label{nmb:4} $C^M$-uniform boundedness principles}
        
        \begin{theorem} (Uniform boundedness principle)\label{nmb:4.1}
                Let $M=(M_k)$ be a non-quasianalytic DC-weight sequence.
                Let $E$, $F$, $G$ be convenient vector spaces and let $U\subseteq F$ be
                $c^\infty$-open. A linear mapping $T:E\to C^M(U,G)$ is bounded if and only
                if $\on{ev}_x\circ T: E\to G$ is bounded for every $x\in U$. 
        \end{theorem}
        
        This is the $C^M$-analogon of \thetag{\ref{nmb:7.3}.7}. Compare with
        \cite[5.22--5.26]{KM97} for the principles behind it. 
        They will be used in the following proof and in \thetag{\ref{nmb:4.6}} and
        \thetag{\ref{nmb:4.10}}. 
        
        \begin{demo}{Proof}
                For $x\in U$ and $\ell\in G^*$ the linear mapping $\ell\circ\on{ev}_x = C^M(x,\ell):C^M(U,G)\to
                \mathbb R$ is continuous, thus $\on{ev}_x$ is bounded. 
                So if $T$ is bounded then so is $\on{ev}_x\circ T$.
                
                Conversely, suppose that $\on{ev}_x\circ T$ is bounded for all $x\in U$. 
                For each closed absolutely convex bounded $B\subseteq E$ we consider the Banach
                space $E_B$. 
                For each $\ell\in
                G^*$, each $C^M$-curve $c:\mathbb R\to U$, each $t\in\mathbb R$, and each compact $K\subset \mathbb R$ 
                the composite given by the following diagram
                is bounded. 
                \[
                \xymatrix{
                        E\ar@{->}[0,1]^{T\qquad} &C^M(U,G)\ar@{->}[1,0]^{C^M(c,\ell)} \ar@{->}[0,2]
                        ^{\on{ev}_{c(t)}} & &G\ar@{->}[1,0]^{\ell} \\
                        E_B\ar@{->}[-1,0] \ar@{->}[0,1] &C^M(\mathbb R,\mathbb R)\ar@{->}[0,1] &
                        \varinjlim _{\rho }C^M_\rho (K,\mathbb R)\ar@{->}[0,1]^{\qquad\on{ev}_t} &
                        \mathbb R \\
                }
                \]
                By \cite[5.24 and 5.25]{KM97} the map $T$ is bounded. In more detail:
                Since $\varinjlim_{\rho}C^M_\rho(K,\mathbb R)$ is webbed by \thetag{\ref{nmb:2.8}}, the
                closed graph theorem \cite[52.10]{KM97} yields that the mapping 
                $E_B\to \varinjlim_{\rho}C^M_\rho(K,\mathbb R)$ is continuous. Thus $T$ is
                bounded.
                \qed\end{demo}
        
        \begin{corollary}\label{nmb:4.2} 
                Let $M=(M_k)$ be a non-quasianalytic DC-weight sequence.
                \begin{enumerate}
                \item
                        For convenient vector spaces $E$ and $F$,
                        on $L(E,F)$ the following bornologies coincide which are induced by:
                        \begin{itemize}
                        \item The topology of uniform convergence on bounded subsets
                                of $E$.
                        \item The topology of pointwise convergence.
                        \item The embedding $L(E,F)\subset C^\infty(E,F)$.
                        \item The embedding $L(E,F)\subset C^M(E,F)$.
                        \end{itemize}
                \item Let $E$, $F$, $G$ be convenient vector spaces and let $U\subset E$ be
                        $c^\infty$-open. A mapping $f:U\times F\to G$ which is linear in
                        the second variable is $C^M$ if and only if 
                        $f^\vee:U\to L(F,G)$ is well defined and $C^M$.
                \end{enumerate}
                Analogous results hold for spaces of multilinear mappings. 
        \end{corollary}
        
        \demo{Proof}
        {\rm (1)}
        That the first three topologies on $L(E,F)$ have the same bounded sets has been shown
        in \cite[5.3 and 5.18]{KM97}.
        The inclusion $C^M(E,F)\to C^\infty(E,F)$ is bounded by \thetag{\ref{nmb:3.10}}
        and by the uniform boundedness principle in \thetag{\ref{nmb:7.3}.7}.
        It remains to show that
        the inclusion $L(E,F)\to C^M(E,F)$ is bounded, where the former space is considered with
        the topology of uniform convergence on bounded sets.
        This follows from the uniform boundedness principle \thetag{\ref{nmb:4.1}}.
        
        {\rm (2)} The assertion for $C^\infty$ is true by \thetag{\ref{nmb:7.3}.6}.
        
        If $f$ is $C^M$ 
        let $c:\mathbb R\to U$ be a $C^M$-curve. We have to show 
        that $t\mapsto f^\vee(c(t))\in L(F,G)$ is $C^M$. By \cite[5.18]{KM97} and
        \thetag{\ref{nmb:3.5}} it suffices to show that 
        $t\mapsto \ell(f^\vee(c(t))(v))=\ell(f(c(t),v))\in \mathbb R$ is $C^M$ for each
        $\ell\in G^*$ and $v\in F$; this is obviously true. 
        
        Conversely, let $f^\vee:U\to L(F,G)$ be $C^M$. 
        We claim that $f:U\times F\to G$ is $C^M$. By composing with $\ell\in G^*$ we
        may assume that $G=\mathbb R$.
        By induction we have 
        \begin{align*}
                &d^kf(x,w_0)\big((v_k,w_k),\dots,(v_1,w_1)\big) 
                = d^k (f^\vee)(x)(v_k,\dots,v_1)(w_0) 
                +\\&
                + \sum_{i=1}^k d^{k-1}(f^\vee)(x)(v_k,\dots,\widehat{v_i,}\dots,v_1)(w_i)
        \end{align*}
        We check condition \thetag{\ref{nmb:3.9}.3} for $f$:
        \begin{align*}
                &\|d^k f (x,w_0)\|_{L^k(E_B\times F_{B'},\mathbb R)}\le
                \\&
                \le \|d^k (f^\vee)(x)(\dots)(w_0)\|_{L^k(E_B,\mathbb R)}
                + \sum_{i=1}^k \|d^{k-1}(f^\vee)(x)\|_{L^{k-1}(E_B,L(F_{B'},\mathbb R))}
                \\&
                \le \|d^k (f^\vee)(x)\|_{L^k(E_B,L(F_{B'},\mathbb R))}\|w_0\|_{B'}
                + \sum_{i=1}^k \|d^{k-1}(f^\vee)(x)\|_{L^{k-1}(E_B,L(F_{B'},\mathbb R))}
                \\&
                \le C\, \rho^k\, k!\, M_k \|w_0\|_{B'}
                + \sum_{i=1}^k C \,\rho^{k-1}\, (k-1)!\, M_{k-1} 
                = C\,\rho^k \,k!\,M_k (\|w_0\|_{B'}+ \tfrac{M_{k-1}}{\rho\, M_k})
        \end{align*}
        where we used \thetag{\ref{nmb:3.9}.3} for $L(i_{B'},\mathbb R)\circ f^\vee: U\to
        L(F_{B'},\mathbb R)$.
        Thus $f$ is $C^M$.
        \qed\enddemo
        
        \begin{proposition}\label{nmb:4.3}
                Let $M=(M_k)$ be a non-quasianalytic DC-weight sequence.
                Let $E$ and $F$ be convenient vector spaces and let $U\subseteq E$ be $c^\infty$-open.
                Then we have the bornological identity
                \[
                        C^M(U,F)=\varprojlim_{s} C^M(\mathbb R,F),
                \]
                where $s$ runs through the strongly $C^M$-curves in $U$ and the connecting mappings are given by 
                $g^*$ for all reparametrizations $g\in C^M(\mathbb R,\mathbb R)$ of curves $s$.
        \end{proposition}
        
        \begin{demo}{Proof}
                By \thetag{\ref{nmb:3.9}} the linear spaces $C^M(U,F)$, $\varprojlim_{s} C^M(\mathbb R,F)$ and 
                $\varprojlim_{c} C^M(\mathbb R,F)$ coincide, where $c$ runs through the
                $C^M$-curves in $U$: 
                Each element $(f_c)_c$ determines a unique function $f:U\to F$ given by 
                $f(x):=(f\circ\on{const}_x)(0)$ with $f\circ c=f_c$ for all such curves $c$, and 
                $f\in C^M$ if and only if $f_c\in C^M$ for all such $c$, by \thetag{\ref{nmb:3.9}}.
                
                Since $C^M(\mathbb R,F)$ carries the initial structure with respect to $\ell_*$ for all $\ell\in F^*$
                we may assume $F=\mathbb R$.
                Obviously the identity $\varprojlim_{c} C^M(\mathbb R,\mathbb R) \to \varprojlim_{s} C^M(\mathbb R,\mathbb R)$ is 
                continuous.
                As projective limit the later space is convenient,
                so we may apply the uniform boundedness principle \thetag{\ref{nmb:4.1}} to conclude that
                the identity in the converse direction is bounded.
                \qed\end{demo}
        
        \begin{proposition}\label{nmb:4.4}
                Let $M=(M_k)$ be a non-quasianalytic DC-weight sequence.
                Let $E$ and $F$ be convenient vector spaces and let $U\subseteq E$ be $c^\infty$-open.
                Then the bornology of $C^M(U,F)$ is initial with respect to each of the following families of mappings
                \begin{align*}
                        i_B^*=C^M(i_B,F)&:C^M(U,F)\to C^M(U\cap E_B,F) \tag{1},\\
                        C^M(i_B,\pi_V)&:C^M(U,F)\to C^M(U\cap E_B,F_V) \tag{2},\\
                        C^M(i_B,\ell)&:C^M(U,F)\to C^M(U\cap E_B,\mathbb R) \tag{3},
                \end{align*}
                where $B$ runs through the closed absolutely convex bounded subsets of $E$
                and $i_B:E_B\to E$ denotes the inclusion, and where $\ell$ runs through the
                continuous linear functionals on $F$, 
                and where $V$ runs through the absolutely convex 0-neighborhoods of $F$ and
                $F_V$ is obtained by factoring out the kernel of the Minkowsky functional of $V$ 
                and then taking the completion with respect to the induced norm.
        \end{proposition}
        
        \noindent{\bf Warning:}
        The structure in {\rm (2)} gives a projective limit description of 
        $C^M(U,F)$ if and only if $F$ is complete since then $F=\varprojlim_V F_V$.
        
        \begin{demo}{Proof}
                Since $i_B:E_B\to E$, $\pi_V:F\to F_V$ and $\ell:F\to\mathbb R$ are bounded linear
                the mappings $i_B^*$, $C^M(i_B,\pi_V)$ and $C^M(i_B,\ell)$ 
                are bounded and linear.
                
                The structures given by {\rm (1)}, {\rm (2)} and {\rm (3)} are 
                successively weaker.
                So let, conversely, $C^M(i_B,\ell)(B)$ be bounded in $C^M(U\cap E_B,\mathbb R)$ for all $B$ 
                and $\ell$.
                By \thetag{\ref{nmb:4.3}} $C^M(U,F)$ carries the initial structure with respect to all 
                $c^*:C^M(U,F)\to C^M(\mathbb R,F)$, where $c:\mathbb R\to U$ are the strongly $C^M$ curves
                and these factor locally as (strongly) $C^M$-curves into some $E_B$.
                By definition $C^M(\mathbb R,F)$ carries the initial structure with respect to 
                $C^M(\iota_I,\ell):C^M(\mathbb R,F)\to C^M(I,\mathbb R)$ where $\iota_I:I\hookrightarrow\mathbb R$ are the inclusions 
                of compact intervals into $\mathbb R$ and $\ell\in F^*$.
                Thus $C^M(U,F)$ carries the initial structure with respect to 
                $C^M(c|_I,\ell):C^M(U,F)\to C^M(I,\mathbb R)$,
                which is coarser than that induced by $C^M(U,F)\to C^M(U\cap E_B,\mathbb R)$.
                \qed\end{demo}
        
        \subsection{\label{nmb:4.5}Definition}
        Let $E$ and $F$ be Banach spaces
        and $A\subseteq E$ convex.
        We consider the linear space
        $C^\infty(A,F)$ consisting of all sequences 
        $(f^k)_{k}\in\prod_{k\in\mathbb N} C(A,L^k(E,F))$ satisfying
        \[f^k(y)(v)-f^k(x)(v) =
        \int_0^1 f^{k+1}(x+t(y-x))(y-x,v)\,dt 
        \]
        for all $k\in\mathbb N$, $x,y\in A$, and $v\in E^k$.
        If $A$ is open we can identify this space with that of all smooth functions
        $A\to F$ by passing to jets.
        
        In addition, let $M=(M_k)$ be a non-quasianalytic DC-weight sequence and $(r_k)$ a sequence of 
        positive real numbers. Then we consider the normed spaces 
        \[
                C^M_{(r_k)}(A,F):=\Bigl\{(f^k)_k\in C^\infty(A,F):\|(f^k)\|_{(r_k)}<\infty\Bigr\}
        \]
        where the norm is given by
        \[
                \|(f^k)\|_{(r_k)}:=\sup\Bigl\{\frac{\|f^k(a)(v_1,\dots,v_k)\|}{k!\,r_k\,M_k\,\|v_1\|\cdot\dots\cdot \|v_k\|}:
                k\in\mathbb N, a\in A, v_i\in E\Bigr\}.
        \]
        If $(r_k)=(\rho^k)$ for some $\rho>0$ we just write $\rho$ instead of $(r_k)$ as indices.
        The spaces $C^M_{(r_k)}(A,F)$ are Banach spaces, since they are closed
        in $\ell^\infty(\mathbb N,\ell^\infty(A,L^k(E,F)))$ via 
        $(f^k)_k\mapsto (k\mapsto \frac1{k!\,r_k\,M_k} f^k)$.
        
        \begin{theorem}\label{nmb:4.6} 
                Let $M=(M_k)$ be a non-quasianalytic DC-weight sequence.
                Let $E$ and $F$ be Banach spaces and let $U\subseteq E$ be open.
                Then the space $C^M(U,F)$ can be described bornologically 
                in the following equivalent ways, i.e.\ these constructions give the same
                vector space and the same bounded sets.
                \begin{align*}
                        \varprojlim_K \varinjlim_{\rho,W} C^M_\rho(W,F) \tag{1} \\
                        \varprojlim_K \varinjlim_\rho C^M_\rho(K,F) \tag{2} \\
                        \varprojlim_{K,(r_k)} C^M_{(r_k)}(K,F) \tag{3} \\
                        \varprojlim_{c,I} \varinjlim_\rho C^M_\rho(I,F) \tag{4}
                \end{align*}
                Moreover, all involved inductive limits are regular, i.e.\ the bounded sets of the inductive limits 
                are contained and bounded in some step.
                
                Here $K$ runs through all compact convex subsets of $U$ ordered by inclusion, 
                $W$ runs through the open subsets $K\subseteq W\subseteq U$ again ordered by inclusion,
                $\rho$ runs through the positive real numbers,
                $(r_k)$ runs through all sequences of positive real numbers for which $\rho^k/r_k\to 0$ for all 
                $\rho>0$,
                $c$ runs through the $C^M$-curves in $U$ ordered by reparametrization with $g\in C^M(\mathbb R,\mathbb R)$
                and $I$ runs through the compact intervals in $\mathbb R$.
        \end{theorem}
        
        \begin{demo}{Proof}
                Note first that all four descriptions describe 
                smooth functions $f:U\to F$, which are given by 
                $x\mapsto f^0(x)$ in {\rm (1)}--{\rm (3)} for appropriately chosen $K$ with 
                $x\in K$ where $f^0:K\to F$
                and by
                $x\mapsto f_c(t)$ in {\rm (4)} for $c$ with $x=c(t)$, $t\in I$ and $f_c:I\to F$.
                Smoothness of $f$ follows, since we may test with $C^M$-curves and these factor locally into some 
                $K$.
                
                By \thetag{\ref{nmb:3.9}} all four descriptions describe $C^M(U,F)$ as vector space.

                Obviously the identity is continuous from {\rm (1)} to {\rm (2)} and from 
                {\rm (2)} to {\rm (3)}.
                
                The identity from {\rm (3)} to {\rm (1)} is continuous, since 
                the space given by {\rm (3)} is as inverse limit of Banach spaces convenient and 
                the inductive limit 
                in {\rm (1)} is by construction an (LB)-space, hence webbed, and thus 
                we can apply the uniform $\mathcal{S}$-boundedness principle
                \cite[5.24]{KM97},
                where $\mathcal{S}=\{\on{ev}_x:x\in U\}$.
                
                So the descriptions in {\rm (1)}--{\rm (3)} describe the same
                complete bornology on $C^M(U,F)$ 
                and satisfy the uniform $\mathcal{S}$-boundedness principle.
                
                Moreover, the inductive limits involved in {\rm (1)} and
                {\rm (2)} are regular:
                In fact the bounded sets $\mathcal{B}$ therein are also bounded in the
                structure of {\rm (3)}, i.e.,
                for every compact $K\subseteq U$ and sequence $(r_k)$ of positive real numbers
                for which $\rho^k/r_k\to 0$ for all
                $\rho>0$:
                \[
                        \sup\Bigl\{\frac{\|f^k(a)(v_1,\dots,v_k)\|}{k!\,r_k\,M_k\,\|v_1\|\cdot\dots\cdot \|v_k\|}:
                        k\in\mathbb N, a\in A, v_i\in E, f\in\mathcal{B}\Bigr\} <\infty
                \]
                and so the sequence 
                \[
                        a_k:=
                        \sup\Bigl\{\frac{\|f^k(a)(v_1,\dots,v_k)\|}{k!\,M_k\,\|v_1\|\cdot\dots\cdot \|v_k\|}:
                        a\in A, v_i\in E, f\in\mathcal{B}\Bigr\} <\infty
                \]
                satisfies $\sup_k a_k/r_k<\infty$ for all $(r_k)$ as above. By \cite[9.2]{KM97}
                these are the coefficients of a power series with positive 
                radius of convergence.
                Thus $a_k/\rho^k$ is bounded for some $\rho>0$. 
                This means that $\mathcal{B}$ is contained and bounded in $C^M_\rho(K,F)$. 
                
                That also {\rm (4)} describes the same bornology follows again by 
                the $\mathcal{S}$-uniform boundedness principle, 
                since the inductive limit in {\rm (4)} is regular by what 
                we said before for the special 
                case $E=\mathbb R$ and hence the structure of {\rm (4)} is convenient.
                \qed\end{demo}

        \begin{lemma}\label{nmb:4.7}
                Let $M$ be a non-quasianalytic DC-weight sequence. 
                For any convenient
                vector space $E$ the flip of variables induces an 
                isomorphism $L(E,C^M(\mathbb R,\mathbb R)) \cong 
                C^M(\mathbb R,E')$ as vector spaces. 
        \end{lemma}
        
        \begin{demo}{Proof} For $c\in C^M(\mathbb R,E')$ consider $\tilde c(x)
                := \on{ev}_x\circ c\in C^M(\mathbb R,\mathbb R)$ for $x\in E$. By the
                uniform boundedness principle \thetag{\ref{nmb:4.1}} 
                the linear mapping $\tilde c$ is bounded, since 
                $\on{ev}_t\circ\tilde c=c(t)\in E'$.
                
                If conversely $\ell\in L(E,C^M(\mathbb R,\mathbb R))$, we consider
                $\tilde \ell(t) = \on{ev}_t\circ \ell \in E'=L(E,\mathbb R)$ for $t\in \mathbb R$. 
                Since the bornology of $E'$ is generated by $\mathcal{S}:=\{ev_x:x\in E\}$, 
                $\tilde \ell:\mathbb R \to E'$ is $C^M$, for 
                $\on{ev}_x\circ\tilde \ell=\ell(x)\in C^M(\mathbb R,\mathbb R)$, by \thetag{\ref{nmb:3.5}}.
                \qed\end{demo}
        
        \begin{lemma}\label{nmb:4.8}
                Let $M=(M_k)$ be a non-quasianalytic DC-weight sequence.
                By $\lambda^M(\mathbb R)$ we denote the $c^\infty$-closure of the linear subspace
                generated by $\{\on{ev}_t:t\in\mathbb R\}$ in $C^M(\mathbb R,\mathbb R)'$ and let $\delta:\mathbb R\to\lambda^M(\mathbb R)$
                be given by $t\mapsto \on{ev}_t$.
                Then $\lambda^M(\mathbb R)$ is the free convenient vector space over $C^M$, i.e.\
                for every convenient vector space $G$ the $C^M$-curve $\delta$ induces a
                bornological isomorphism
                \begin{align*}
                        L(\lambda^M(\mathbb R),G)&\cong C^M(\mathbb R,G).
                \end{align*}
        \end{lemma}
        We expect $\lambda^M(\mathbb R)$ to be equal to $C^M(\mathbb R,\mathbb R)'$ as it is the case 
        for the analogous situation of smooth mappings, see \cite[23.11]{KM97}, and of
        holomorphic mappings, see \cite{Siegl95} and \cite{Siegl97}.
        \begin{demo}{Proof}
                The proof goes along the same lines as in \cite[23.6]{KM97} and in \cite[5.1.1]{FK88}.
                Note first that 
                $\lambda^M(\mathbb R)$ is a convenient vector space since it is $c^\infty$-closed in the convenient vector space 
                $C^M(\mathbb R,\mathbb R)'$. Moreover, $\delta$ is $C^M$ by \thetag{\ref{nmb:3.5}}, 
                since $\on{ev}_h\circ\delta=h$ for all $h\in C^M(\mathbb R,\mathbb R)$, so $\delta^*:L(\lambda^M(\mathbb R),G)\to C^M(\mathbb R,G)$
                is a well-defined linear mapping. 
                This mapping is injective, since each bounded linear mapping $\lambda^M(\mathbb R)\to G$ is uniquely 
                determined on $\delta(\mathbb R)=\{\on{ev}_t:t\in\mathbb R\}$.
                Let now $f\in C^M(\mathbb R,G)$. Then 
                $\ell\circ f\in C^M(\mathbb R,\mathbb R)$ for every $\ell\in G^*$ and hence 
                $\tilde f:C^M(\mathbb R,\mathbb R)'\to \prod_{G^*}\mathbb R$ given by $\tilde f(\varphi)=(\varphi(\ell\circ f))_{\ell\in G^*}$
                is a well-defined bounded linear map. Since it maps $\on{ev_t}$ to $\tilde f(\on{ev}_t)=\delta(f(t))$, 
                where $\delta:G\to\prod_{G^*}\mathbb R$ denotes the bornological embedding given by 
                $x\mapsto (\ell(x))_{\ell\in G^*}$, it induces a bounded linear mapping $\tilde f:\lambda^M(\mathbb R)\to G$
                satisfying $\tilde f\circ\delta=f$.
                Thus $\delta^*$ is a linear bijection.
                That it is a bornological isomorphism (i.e.\ $\delta^*$ and its inverse are both bounded) follows from 
                the uniform boundedness principles \thetag{\ref{nmb:4.1}} and \thetag{\ref{nmb:4.2}}.
                \qed\end{demo}
        
        \begin{corollary}\label{nmb:4.9}
                Let $M=(M_k)$ and $N=(N_k)$ be non-quasianalytic DC-weight sequences.
                We have the following isomorphisms of linear spaces
                \begin{enumerate}
                \item $C^\infty(\mathbb R,C^M(\mathbb R,\mathbb R)) \cong 
                        C^M(\mathbb R,C^\infty(\mathbb R,\mathbb R))$ 
                \item $C^\omega(\mathbb R,C^M(\mathbb R,\mathbb R)) \cong 
                        C^M(\mathbb R,C^\omega(\mathbb R,\mathbb R))$
                \item $C^N(\mathbb R,C^M(\mathbb R,\mathbb R)) \cong 
                        C^M(\mathbb R,C^N(\mathbb R,\mathbb R))$
                \end{enumerate}
        \end{corollary}
        
        \begin{demo}{Proof} 
                For $\alpha\in\{\infty,\omega,N\}$ we get
                \begin{align*}
                        C^M(\mathbb R,C^\alpha(\mathbb R,\mathbb R)) &\cong
                        L(\lambda^M(\mathbb R),C^\alpha(\mathbb R,\mathbb R))\qquad\text{ by \thetag{\ref{nmb:4.8}}}\\
                        &\cong C^\alpha(\mathbb R,L(\lambda^M(\mathbb R),\mathbb R))\qquad\text{by \thetag{\ref{nmb:4.7}}, 
                                \cite[3.13.4, 5.3, 11.15]{KM97}}\\
                        &\cong C^\alpha(\mathbb R,C^M(\mathbb R,\mathbb R))\qquad\text{ by \thetag{\ref{nmb:4.8}}.}\qed
                \end{align*}
        \end{demo}
        
        \begin{theorem}\label{nmb:4.10} (Canonical isomorphisms)
                Let $M=(M_k)$ and $N=(N_k)$ be non-quasianalytic DC-weight sequences.
                Let $E$, $F$ be convenient vector spaces and let
                $W_i$ be $c^\infty$-open subsets in such. 
                We have the following natural bornological isomorphisms:
                \begin{enumerate}
                \item[(1)] $C^M(W_1,C^N(W_2,F))\cong C^N(W_2,C^M(W_1,F))$,
                \item[(2)] $C^M(W_1,C^\infty(W_2,F))\cong C^\infty(W_2,C^M(W_1,F))$.
                \item[(3)] $C^M(W_1,C^\omega(W_2,F))\cong C^\omega(W_2,C^M(W_1,F))$.
                \item[(4)] $C^M(W_1,L(E,F))\cong L(E,C^M(W_1,F))$.
                \item[(5)] $C^M(W_1,\ell^\infty(X,F))\cong \ell^\infty(X,C^M(W_1,F))$.
                \item[(6)] $C^M(W_1,\operatorname{\mathcal Lip}^k(X,F))\cong \operatorname{\mathcal Lip}^k(X,C^M(W_1,F))$.
                \end{enumerate}
                In {\rm (5)} the space $X$ is an $\ell^\infty$-space, i.e.\ 
                a set together 
                with a bornology induced by a family of real valued functions on $X$, 
                cf.\ \cite[1.2.4]{FK88}.
                In {\rm (6)} the space $X$ is a $\operatorname{\mathcal Lip}^k$-space, cf.\ 
                \cite[1.4.1]{FK88}.
                The spaces $\ell^\infty(X,F)$ and $\operatorname{\mathcal Lip}^k(W,F)$ are defined in 
                \cite[3.6.1 and 4.4.1]{FK88}. 
        \end{theorem}
        
        \begin{demo}{Proof} 
                All isomorphisms, as well as their inverse mappings, are 
                given by the flip of coordinates: $f\mapsto \tilde f$, where 
                $\tilde f(x)(y):=f(y)(x)$. Furthermore, all occurring 
                function spaces are convenient and satisfy the uniform 
                $\mathcal{S}$-boundedness theorem, where $\mathcal{S}$ is the set of point 
                evaluations, by \thetag{\ref{nmb:4.1}}, \cite[11.11, 11.14, 11.12]{KM97}, 
                and by 
                \cite[3.6.1, 4.4.2, 3.6.6, and 4.4.7]{FK88}.
                
                That $\tilde f$ has values in the corresponding spaces follows 
                from the equation $\tilde f(x)=ev_x\circ f$. 
                One only has to check that $\tilde f$ itself 
                is of the corresponding class, since it follows that $f\mapsto \tilde 
                f$ is bounded. This is a consequence of the uniform boundedness 
                principle, since 
                \begin{displaymath}
                        (\on{ev}_x\circ \tilde{(\quad)})(f) = \on{ev}_x(\tilde f)= 
                        \tilde f(x) = \on{ev}_x\circ f = (\on{ev}_x)_*(f).\end{displaymath}
                
                That $\tilde f$ is of the appropriate class in {\rm (1)} %
                and in {\rm (2)} follows by composing with the appropriate curves
                $c_1:\mathbb R\to W_1$, $c_2:\mathbb R\to W_2$ and $\lambda\in F^*$
                and thereby reducing the statement to the special case in \thetag{\ref{nmb:4.9}}.
                
                That $\tilde f$ is of the appropriate class in {\rm (3)} follows
                by composing with
                $c_1 \in C^M(\mathbb R,W_1)$ and 
                $C^{\beta_2}(c_2,\lambda):C^\omega(W_2,F) \to C^{\beta_2}(\mathbb R,\mathbb R)$ for 
                all $\lambda\in F^*$ and $c_2\in C^{\beta_2}(\mathbb R,W_2)$, where $\beta_2$ 
                is in $\{\infty,\omega\}$.
                Then 
                $C^{\beta_2}(c_2,\lambda)\circ \tilde f \circ c_1 =
                (C^M(c_1,\lambda)\circ f \circ c_2)^\sim
                : \mathbb R \to C^{\beta_2}(\mathbb R,\mathbb R)$
                is $C^M$ by \thetag{\ref{nmb:4.9}}, since 
                $C^M(c_1,\lambda)\circ f \circ c_2
                : \mathbb R \to W_2 \to C^M(W_1,F) \to C^M(\mathbb R,\mathbb R)$
                is $C^{\beta_2}$.
                
                That $\tilde f$ is of the appropriate class in {\rm (4)} follows,
                since $L(E,F)$ is the $c^\infty$-closed subspace of $C^M(E,F)$ formed
                by the linear $C^M$-mappings.
                
                That $\tilde f$ is of the appropriate class in {\rm (5)} or 
                {\rm (6)} follows from
                {\rm (4)}, using the free convenient vector spaces $\ell^1(X)$ 
                or $\lambda^k(X)$ over 
                the $\ell^\infty$-space $X$ or the the $\operatorname{\mathcal Lip}^k$-space $X$, see 
                \cite[5.1.24 or 5.2.3]{FK88}, 
                satisfying $\ell^\infty(X,F)\cong L(\ell^1(X),F)$ or 
                satisfying $\operatorname{\mathcal Lip}^k(X,F)\cong L(\lambda^k(X),F)$. 
                Existence of these free convenient vector spaces can be proved in a 
                similar way as in \thetag{\ref{nmb:4.8}}.
                \qed
        \end{demo}
        
        \section{\label{nmb:5} Exponential law}
        
        \subsection{\label{nmb:5.1}Difference quotients}
        For the following see \cite[1.3]{FK88}.
        For a subset $K\subseteq \mathbb R^n$, $\alpha=(\alpha_1,\dots,\alpha_n)\in \mathbb{N}^n$, 
        a linear space $E$, and $f:K\to E$ let:
        \begin{align*}
                &\mathbb R^{\langle k \rangle} = \bigl\{(x_0,\dots,x_k)\in \mathbb R^{k+1}:
                x_i\ne x_j\text{ for }i\ne j\bigr\}
                \\
                &K^\alpha=\bigl\{(x^1,\dots,x^n)\in \mathbb R^{\alpha_1+1}\times\dots\times\mathbb
                R^{\alpha_n+1}: 
                (x^1_{i_1},\dots,x^n_{i_n})\in K\text{ for }0\le i_j\le \alpha_j\bigr\}
                \\
                &K^{\langle \alpha \rangle} = K^{\alpha}\cap (\mathbb R^{\langle \alpha_1
                        \rangle}\times\dots\times\mathbb R^{\langle \alpha_n\rangle})
                \\
                &\beta_i(x) = k! \prod_{\substack{ 0\le j\le k \\ j\ne i}} \frac1{x_i-x_j}
                \text{ for }x=(x_0,\dots,x_k)\in \mathbb R^{\langle k \rangle} 
                \\
                &\delta^{\alpha}f(x^1,\dots,x^n) = \sum_{i_1=0}^{\alpha_1}\dots\sum_{i_n=0}^{\alpha_n}
                \beta_{i_1}(x^1)\dots\beta_{i_n}(x^n) f(x_{i_1}^1,\dots, x_{i_n}^n)
        \end{align*}
        Note that $\delta^0f=f$ and $\delta^\alpha = \delta_n^{\alpha_n}\circ\dots\circ\delta_1^{\alpha_1} $where
        \[
        \delta_i^k g(x^1,\dots,x^n) =
        \delta^k(g(x^1,\dots,x^{i-1},\quad,x^{i+1},\dots,x^n))(x^i).
        \]

        \begin{lemma*}
                Let $E$ be a convenient vector space. Let $U\subseteq \mathbb R^n$ be open.
                For $f:U\to E$ 
                the following conditions are equivalent:
                \begin{enumerate}
                \item $f:U\to E$ is $C^M$.
                \item For every compact convex set $K$ in $U$ and every $\ell\in E^*$ there exists 
                        $\rho>0$ such that 
                        \[
                        \left\{\frac{\delta^\alpha (\ell\circ f)(x)}{\rho^{|\alpha|}\, |\alpha|!\, M_{|\alpha|}}:
                        \alpha\in\mathbb N^n, x\in K^{\langle \alpha \rangle}
                        \right\}
                        \]
                        is bounded in $\mathbb R$.
                \end{enumerate}
                Furthermore, the norm on the space $C^M_\rho(K,\mathbb R)$ from \thetag{\ref{nmb:2.8}} (for
                convex $K$) is also given by
                \[
                \|f\|_{\rho,K} :=
                \sup \Big\{ \frac{|\delta^{\alpha} f(x)|}{\rho^{|\alpha|} |\alpha|! M_{|\alpha|}} : \alpha \in
                \mathbb{N}^n, x \in K^{\langle \alpha\rangle}\Big\}.
                \]
        \end{lemma*}
        
        \begin{demo}{Proof}
                By composing with bounded linear functionals we may assume that $E=\mathbb
                R$.
                
                {\rm (1)} $\implies$ {\rm (2)}
                If $f$ is $C^M$ then for each compact convex set $K$ in $U$
                there exists $\rho>0$ such that 
                \[
                \left\{\tfrac{\partial^\alpha f (x)}{\rho^{|\alpha|}\, {|\alpha|}! M_{|\alpha|}}: 
                \alpha\in\mathbb N^n, x\in K \right\}
                \]
                is bounded in $\mathbb R$. 
                
                For a differentiable function $g:\mathbb R\to \mathbb R$ and $t_0<\dots<t_j$ there 
                exist $s_i$ with $t_i<s_i<t_{i+1}$ such that
                \[
                \delta^{j}g(t_0,\dots,t_j)=\delta^{j-1}g'(s_0,\dots,s_{j-1}).
                \]
                This follows by Rolle's theorem, see \cite[12.4]{KM97}. 
                Recursion, for $g= \partial^\alpha f$, shows that 
                $\delta^{\alpha} f(x^0,\dots,x^n) =\partial^\alpha f(s)$ for some $s\in K$.
                
                {\rm (2)} $\implies$ {\rm (1)}
                $f$ is $C^\infty$ by \cite[1.3.29]{FK88} since 
                each difference quotient $\delta^{\alpha} f$ is bounded on bounded sets.
                
                For $g\in C^\infty(\mathbb R,\mathbb R)$, using (see \cite[1.3.6]{FK88})
                \[
                g(t_j) = \sum_{i=0}^j\tfrac1{i!}
                \prod_{l=0}^{i-1}(t_j-t_l)\;\delta^jg(t_0,\dots,t_j),
                \]
                induction on $j$ and differentiability of $g$ shows that 
                \[\delta^jg'(t_0,\dots,t_j)
                =\tfrac1{j+1}\sum_{i=0}^j\delta^{j+1}g(t_0,\dots,t_j,t_i),
                \]
                where 
                $\delta^{j+1}g(t_0,\dots,t_j,t_i)
                := \lim_{t\to t_i}\delta^{j+1}g(t_0,\dots,t_j,t)$.
                If the right hand side divided by $\rho^{|\alpha|}\,|\alpha|!\,M_{|\alpha|}$ is bounded, 
                then also $\delta^jg'/(\rho^{|\alpha|}\,|\alpha|!\,M_{|\alpha|})$ is bounded. 
                
                By recursion, applied to $g=\delta^\beta\partial^{\alpha-\beta} f$, we conclude that $f\in C^M$.
                \qed\end{demo}
        
        \begin{lemma}\label{nmb:5.2}
                Let $E$ be a convenient vector space such that there exists a Baire vector
                space topology on the dual $E^*$ for which the point evaluations $\on{ev}_x$
                are continuous for all $x\in E$. For a mapping $f:\mathbb R^n\to E$ the
                following are equivalent:
                \begin{enumerate}
                \item $\ell\circ f$ is $C^M$ for all $\ell\in E^*$.
                \item For every convex compact $K\subseteq \mathbb R^n$ there exists
                        $\rho>0$ such that 
                        \[
                        \left\{\frac{\partial^\alpha f(x)}{\rho^{|\alpha|}\, |\alpha|!\, M_{|\alpha|}}: \alpha\in\mathbb{N}^n, x\in K
                        \right\}\text{ is bounded in }E.
                        \]
                \item For every convex compact $K\subseteq \mathbb R^n$ there exists 
                        $\rho>0$ such that 
                        \[
                        \left\{\frac{\delta^\alpha f(x)}{\rho^{|\alpha|}\, |\alpha|!\, M_{|\alpha|}}: \alpha\in\mathbb{N}^n, x\in
                        K^{\langle \alpha\rangle} \right\}\text{ is bounded in }E.
                        \]
                \end{enumerate}
        \end{lemma}
        
        \begin{demo}{Proof}
                {\rm (2)} $\implies$ {\rm (1)} is obvious.
                
                {\rm (1)} $\implies$ {\rm (2)}
                Let $K$ be compact convex in
                $\mathbb R^n$. 
                We consider the sets
                \[
                A_{\rho,C} :=\Bigl\{\ell\in E^*: \frac{|\partial^\alpha (\ell\circ f)(x)|}{\rho^{|\alpha|}\,
                        |\alpha|!\, M_{|\alpha|}}\le C\text{ for all }\alpha\in \mathbb{N}^n, x\in K\Bigr\}
                \]
                which are closed subsets in $E^*$ for the Baire topology. We have
                $\bigcup_{\rho,C}A_{\rho,C}= E^*$. By the Baire property there exists $\rho$
                and $C$ such that the interior $U$ of $A_{\rho,C}$ is non-empty. If
                $\ell_0\in U$ then for all $\ell\in E^*$ there is an $\epsilon>0$ such that 
                $\epsilon\ell\in U-\ell_0$ and hence for all $x\in K$ and all $\alpha$ we have
                \begin{align*}
                        |\partial^\alpha (\ell\circ f)(x)| \le \tfrac1\epsilon \left(|\partial^\alpha ((\epsilon\ell+\ell_0)\circ f)(x)| +
                        |\partial^\alpha (\ell_0\circ f)(x)|\right) \le \tfrac{2C}{\epsilon}\,\rho^{|\alpha|}\,|\alpha|!\,
                        M_{|\alpha|}.
                \end{align*}
                So the set 
                \[
                \left\{\frac{\partial^\alpha f(x)}{\rho^{|\alpha|}\, |\alpha|!\, M_{|\alpha|}}: \alpha\in\mathbb{N}^n, x\in K
                \right\}
                \]
                is weakly bounded in $E$ and hence bounded. 
                
                {\rm (3)} $\implies$ {\rm (1)} follows by Lemma
                \thetag{\ref{nmb:5.1}}. 
                {\rm (1)} $\implies$ {\rm (3)} follows as above for the
                difference quotients instead of the partial differentials. 
                \qed\end{demo}

        \begin{theorem}\label{nmb:5.3}
                (Cartesian closedness) 
                Let $M=(M_k)$ be a non-quasianalytic DC-weight sequence of moderate growth
                \eqref{mgrowth}.
                Then the
                category of $C^M$-mappings between convenient real vector spaces 
                is cartesian closed.
                More precisely, for convenient vector spaces $E$, $F$ and $G$
                and $c^\infty$-open sets $U\subseteq E$ and $W\subseteq F$ 
                a mapping $f:U\times W\to G$ is $C^M$ if and only
                if $f^\vee :U\to C^M(W,G)$ is $C^M$.
        \end{theorem}
        
        \begin{demo}{Proof}
                We first show the result for $U=\mathbb R$, $W=\mathbb R$, $G=\mathbb R$.
                
                If $f\in C^M(\mathbb R^2,\mathbb R)$ then clearly for any 
                $x\in \mathbb R$ the function
                $f^\vee (x)=f(x,\quad)\in C^M(\mathbb R,\mathbb R)$.
                To show that $f^\vee :\mathbb R\to C^M(\mathbb R,\mathbb R)$ is $C^M$
                it suffices to check \thetag{\ref{nmb:5.1}.2} for all $\ell\in
                C^M(\mathbb R,\mathbb R)^*$. Such an $\ell$ factors over 
                $\varinjlim_\rho C_\rho^M(L)$ for some compact $L\subset\mathbb R$.
                Let $K\subset \mathbb R$ be compact.
                Since $f$ is $C^M$ there exists $C>0$ and $\rho>0$ by lemma \thetag{\ref{nmb:5.1}} 
                such that 
                \[
                \frac{|\delta^\alpha f(x,y)|}{\rho^{|\alpha|}|\alpha|!M_{|\alpha|}}\leq C\quad\text{ for }
                \alpha\in\mathbb{N}^2, (x,y)\in (K\times L)^{\langle \alpha \rangle}.
                \]
                Since $M$ is of moderate growth \eqref{mgrowth} we have
                $M_{j+k}\le \sigma^{j+k} M_j M_k$ for some $\sigma>0$. 
                Let $\alpha=(\alpha_1,\alpha_2)\in \mathbb{N}^2$. Then:
                \begin{align*}
                        &\left\| \frac{\delta^{\alpha_1}f^\vee (x)}{\rho_1^{\alpha_1}\,\alpha_1!\,M_{\alpha_1}}\right\|_{\rho_2,L}
                        =
                        \sup\Bigl\{\frac{|\delta_2^{\alpha_2}\delta_1^{\alpha_1}f(x,y)|}{\rho_1^{\alpha_1}\,\alpha_1!\,M_{\alpha_1}\,
                                \rho_2^{\alpha_2}\,\alpha_2!\,M_{\alpha_2}}: \alpha_2\in\mathbb{N}, y\in L^{\langle \alpha_2 \rangle}
                        \Bigr\}
                        \\&
                        \le 
                        \sup\Bigl\{\frac{|\delta_2^{\alpha_2}\delta_1^{\alpha_1}f(x,y)|}
                        {\rho_1^{\alpha_1}\,\rho_2^{\alpha_2}
                                \frac{\alpha_1!\,\alpha_2!}{(\alpha_1+\alpha_2)!}\,(\alpha_1+\alpha_2)!
                                \,\sigma^{-\alpha_1-\alpha_2} M_{\alpha_1+\alpha_2}}: \alpha_2\in\mathbb{N}, y\in L^{\langle \alpha_2 \rangle}\Bigr\}
                        \\&
                        \le 
                        \sup\Bigl\{\frac{|\delta^{\alpha}f(x,y)|}
                        {\rho_1^{\alpha_1}\,\rho_2^{\alpha_2}\,\sigma^{-|\alpha|}
                                2^{-|\alpha|}\,|\alpha|!\,
                                M_{|\alpha|}}: \alpha_2\in\mathbb{N}, y\in L^{\langle \alpha_2 \rangle} \Bigr\}
                        \\&
                        \le 
                        \sup\Bigl\{\frac{|\delta^{\alpha}f(x,y)|}
                        {\rho^{|\alpha|}\,
                                |\alpha|!\,M_{|\alpha|}}: \alpha_2\in\mathbb{N}, y\in L^{\langle \alpha_2 \rangle} \Bigr\}
                        \le C \text{ for }\alpha_1\in\mathbb{N}, x\in K^{\langle \alpha_1 \rangle}
                \end{align*}
                for $\rho_1=\rho_2=2\sigma\rho$.
                So $f^\vee :K\to C^M_{\rho_2}(L,\mathbb R)$ is $C^M$. Thus $\ell\circ f^\vee $
                is $C^M$.
                
                Conversely, let $f^\vee :\mathbb R\to C^M(\mathbb R,\mathbb R)$ be $C^M$.
                Then $f^\vee :\mathbb R\to \varinjlim_{\rho_2} C^M_{\rho_2}(L,\mathbb R)$ is $C^M$
                for all compact subsets $L\subset \mathbb R$. The dual space
                $(\varinjlim_{\rho_2} C^M_{\rho_2}(L,\mathbb R))^*$ can be equipped with the Baire
                topology of the countable limit $\varprojlim_{\rho_2} C^M_{\rho_2}(L,\mathbb R)^*$
                of Banach spaces. 
                \[
                \xymatrix{
                        \mathbb R \ar[r]^{f^\vee \quad} & C^M(\mathbb R,\mathbb R) \ar[r] &
                        \varinjlim_{\rho_2} C^M_{\rho_2}(L,\mathbb R) \\
                        K \ar@{^(->}[u] \ar[rr]^{f^\vee } & & C^M_{\rho_2}(L,\mathbb R) \ar[u]
                }
                \]
                Thus the mapping $f^\vee :\mathbb R\to \varinjlim_{\rho_2} C^M_{\rho_2}(L,\mathbb
                R)$ is strongly $C^M$ by \thetag{\ref{nmb:5.2}}. Since the inductive limit 
                $\varinjlim_{\rho_2} C^M_{\rho_2}(L,\mathbb R)$ is countable and regular 
                (\cite[7.4 and 7.5]{Floret71} or \cite[52.37]{KM97}), for each compact 
                $K\subset \mathbb R$ there exists $\rho_1>0$ such that 
                the bounded set 
                \[
                \left\{\frac{\partial^{\alpha_1} f^\vee (x)}{\rho_1^{\alpha_1}\, \alpha_1!\,
                        M_{\alpha_1}}: \alpha_1\in\mathbb{N}, x\in K\right\}
                \]
                is contained and bounded in $C^M_{\rho_2}(L,\mathbb R)$ for some
                ${\rho_2}>0$. 
                Thus for $\alpha_1\in\mathbb{N}$ and $x\in K$ we have (using \eqref{logconvex1})
                \begin{align*}
                        \infty &> C := \sup_{\substack{\alpha_1\in\mathbb{N}\\ y\in K}}
                        \left\| \frac{\delta^{\alpha_1}f^\vee (y)}{\rho_1^{\alpha_1}\,\alpha_1!\,M_{\alpha_1}}\right\|_{\rho_2,L}
                        \ge 
                        \left\| \frac{\delta^{\alpha_1}f^\vee (x)}{\rho_1^{\alpha_1}\,\alpha_1!\,M_{\alpha_1}}\right\|_{\rho_2,L}
                        \\&
                        =
                        \sup\Bigl\{\frac{|\delta_2^{\alpha_2}\delta_1^{\alpha_1}f(x,y)|}{\rho_1^{\alpha_1}\,\alpha_1!\,M_{\alpha_1}\,
                                \rho_2^{\alpha_2}\,\alpha_2!\,M_{\alpha_2}}: \alpha_2\in\mathbb{N}, y\in L^{\langle \alpha_2 \rangle} \Bigr\}
                        \\&
                        \ge 
                        \sup\Bigl\{\frac{|\delta_2^{\alpha_2}\delta_1^{\alpha_1}f(x,y)|}
                        {\rho_1^{\alpha_1}\,\rho_2^{\alpha_2}
                                \frac{\alpha_1!\,\alpha_2!}{(\alpha_1+\alpha_2)!}\,(\alpha_1+\alpha_2)!
                                M_{\alpha_1+\alpha_2}}: \alpha_2\in\mathbb{N}, y\in L^{\langle \alpha_2 \rangle} \Bigr\}
                        \\&
                        \ge 
                        \sup\Bigl\{\frac{|\delta^{\alpha}f(x,y)|}
                        {\rho^{|\alpha|}\,|\alpha|!\,
                                M_{|\alpha|}}: \alpha_2\in\mathbb{N}, y\in L^{\langle \alpha_2 \rangle} \Bigr\}
                \end{align*}
                where $\rho =\max(\rho_1,\rho_2)$. Thus $f$ is $C^M$.
                
                Now we consider the general case. 
                Given a $C^M$-mapping $f:U\times W\to G$ we have to show that 
                $f^\vee : U\to C^M(W,G)$ is $C^M$. Any continuous linear functional on
                $C^M(W,G)$ factors over some step mapping 
                $C^M(c_2,\ell):C^M(W,G)\to C^M(\mathbb R,\mathbb R)$ of the cone in
                \thetag{\ref{nmb:3.1}} where $c_2$ is a $C^M$-curve in $W$ and $\ell\in G^*$. 
                So we have to check that 
                $C^M(c_2,\ell)\circ f^\vee \circ c_1: \mathbb R\to C^M(\mathbb R,\mathbb R)$ is $C^M$
                for every $C^M$-curve $c_1$ in $U$. 
                Since 
                $(\ell\circ f\circ (c_1\times c_2))^\vee = C^M(c_2,\ell)\circ f^\vee \circ c_1$
                this follows from the special case proved above.
                
                If $f^\vee:U\to C^M(W,G)$ is $C^M$ then 
                $(\ell\circ f\circ (c_1\times c_2))^\vee = C^M(c_2,\ell)\circ f^\vee \circ c_1$ is $C^M$
                for all $C^M$-curves $c_1:\mathbb R\to U$, $c_2:\mathbb R\to W$ and
                $\ell\in G^*$. By the special case, $f$ is then $C^M$.
                \qed\end{demo}
        
        \subsection{\label{nmb:5.4}Example: Cartesian closedness is wrong in general}
        Let $M$ be a DC-weight sequence which is strongly non-quasianalytic 
        but not of moderate growth.
        For example, $M_k=2^{k^2}$ satisfies this by
        \thetag{\ref{nmb:2.7}}. %
        Then by \thetag{\ref{nmb:2.4}} %
        there exists $f:\mathbb R^2\to \mathbb
        R$ of class $C^M$ with $\partial^\alpha f(0,0) = |\alpha|! \, M_{|\alpha|}$.
        {\it We claim that $f^\vee :\mathbb R\to C^M(\mathbb R,\mathbb R)$ is not
                $C^M$.}
        
        Since $M$ is not of moderate growth there exist $j_n \nearrow \infty$ and $k_n>0$ 
        such that 
        \[
        \left(\frac{M_{k_n+j_n}}{M_{k_n} M_{j_n}}\right)^\frac1{k_n+j_n} \ge n.
        \]
        Consider the linear functional $\ell:C^M(\mathbb R,\mathbb R)\to \mathbb R$
        given by 
        \[
        \ell(g)=\sum_n\frac{g^{(j_n)}(0)}{j_n!\,M_{j_n}\,n^{j_n}}.
        \]
        This functional is continuous since 
        \[
        \left|\sum_n\frac{g^{(j_n)}(0)}{j_n!\,M_{j_n}\,n^{j_n}}\right| 
        \le 
        \sum_n\frac{g^{(j_n)}(0)}{j_n!\, \rho^{j_n}\,M_{j_n}}\frac{\rho^{j_n}}{n^{j_n}}
        \le 
        C(\rho)\,\|g\|_{\rho,[-1,1]} < \infty
        \]
        for suitable $\rho$ where
        \[
        C(\rho):=\sum_n \rho^{j_n} \frac{1}{n^{j_n}} < \infty
        \]
        for all $\rho$. 
        But $\ell\circ f^\vee $is not $C^M$ since
        \begin{align*}
                &\left\|\ell\circ f^\vee \right\|_{\rho_1,[-1,1]} \ge
                \sup_k\frac{1}{\rho_1^k\,k!\,M_k}\sum_n\frac{f^{(j_n,k)}(0,0)}{j_n!\,M_{j_n}n^{j_n}}
                \\&
                \ge
                \sup_n\frac{1}{\rho_1^{k_n}\,k_n!\,M_{k_n}}\frac{f^{(j_n,k_n)}(0,0)}{j_n!\,M_{j_n}\,n^{j_n}}
                \\&
                \ge
                \sup_n\frac{(j_n+k_n)!\,M_{j_n+k_n}}{\rho_1^{k_n}\,k_n!\,j_n!\,M_{k_n}\,M_{j_n}\,n^{j_n}}
                \ge \sup_n \frac{n^{j_n+k_n}}{\rho_1^{k_n}\,n^{j_n}} = \infty
        \end{align*}
        for all $\rho_1>0$.
        
        \begin{theorem} \label{nmb:5.5}
                Let $M$ be a non-quasianalytic DC-weight sequence which is of
                moderate growth. 
                Let $E$, $F$, etc., be convenient vector spaces and let $U$ and $V$ be
                $c^\infty$-open subsets of such. 
                \newline 
                {\rm (1)} The exponential law holds: 
                \[
                C^M(U,C^M(V,G)) \cong C^M(U\times V, G)
                \]
                \indent\indent is a linear $C^M$-diffeomorphism of convenient vector spaces. 
                \newline
                The following canonical mappings are $C^M$.
                \begin{align*}
                        &\operatorname{ev}: C^M(U,F)\times U\to F,\quad 
                        \operatorname{ev}(f,x) = f(x)
                        \tag{2}\\&
                        \operatorname{ins}: E\to C^M(F,E\times F),\quad
                        \operatorname{ins}(x)(y) = (x,y)
                        \tag{3}\\&
                        (\quad)^\wedge :C^M(U,C^M(V,G))\to C^M(U\times V,G)
                        \tag{4}\\&
                        (\quad)^\vee :C^M(U\times V,G)\to C^M(U,C^M(V,G))
                        \tag{5}\\&
                        \operatorname{comp}:C^M(F,G)\times C^M(U,F)\to C^M(U,G)
                        \tag{6}\\&
                        C^M(\quad,\quad):C^M(F,F_1)\times C^M(E_1,E)\to 
                        C^M(C^M(E,F),C^M(E_1,F_1))
                        \tag{7}\\&
                        \qquad (f,g)\mapsto(h\mapsto f\circ h\circ g)
                        \\&
                        \prod:\prod C^M(E_i,F_i)\to C^M(\prod E_i,\prod F_i)
                        \tag{8}\end{align*}
        \end{theorem}
        
        \begin{demo}{Proof}
                {\rm (2)} 
                The mapping associated to $\on{ev}$ via cartesian 
                closedness is the identity on $C^M(U,F)$, which is $C^M$, 
                thus $\on{ev}$ is also $C^M$.
                
                {\rm (3)} 
                The mapping associated to $\on{ins}$ via cartesian
                closedness is the identity on $E\times F$, hence $\on{ins}$ is $C^M$.
                
                {\rm (4)} 
                The mapping associated to $(\quad)^\wedge$ via cartesian closedness 
                is the $C^M$-composition of evaluations 
                $\on{ev}\circ(\on{ev}\times \on{Id}):(f;x,y)\mapsto f(x)(y)$.
                
                {\rm (5)} 
                We apply cartesian closedness twice to get the associated 
                mapping $(f;x;y)\mapsto f(x,y)$, which is just a 
                $C^M$ evaluation mapping.
                
                {\rm (6)} 
                The mapping associated to $\on{comp}$
                via cartesian closedness is $(f,g;x)\mapsto f(g(x))$, which 
                is the $C^M$-mapping $\on{ev}\circ (\on{Id} \times \on{ev})$.
                
                {\rm (7)} 
                The mapping associated to the one in question by applying cartesian 
                closedness twice is $(f,g;h,x)\mapsto g(h(f(x)))$, which is the 
                $C^M$-mapping $\on{ev}\circ(\on{Id}\times \on{ev})\circ(\on{Id}\times \on{Id}\times \on{ev})$.
                
                {\rm (8)} 
                Up to a flip of factors the mapping associated via 
                cartesian closedness is
                the product of the evaluation mappings 
                $C^M(E_i,F_i)\times E_i\to F_i$.
                
                {\rm (1)}
                follows from {\rm (4)} and {\rm (5)}. 
                \qed\end{demo}

        \section{\label{nmb:6} Manifolds of $C^M$-mappings}
        
        \subsection{\label{nmb:6.1}$C^M$-manifolds}
        Let $M=(M_k)$ be a non-quasianalytic DC-weight sequence of moderate growth. 
        A $C^M$-manifold is a smooth manifold such that all chart changings are
        $C^M$-mappings. Likewise for $C^M$-bundles and $C^M$ Lie groups. 
        
        Note that any finite dimensional (always assumed paracompact) 
        $C^\infty$-manifold admits a
        $C^\infty$-diffeomorphic real analytic structure thus also a
        $C^M$-structure. Maybe, any finite dimensional $C^M$-manifold admits a
        $C^M$-diffeomorphic real analytic structure.
        
        \subsection{\label{nmb:6.2}Spaces of $C^M$-sections}
        Let $E\to B$ be a $C^M$ vector bundle (possibly infinite dimensional).
        The space $C^M(B\gets E)$ of all $C^M$ sections is a convenient
        vector space with the structure induced by 
        \begin{gather*}
                C^M(B\gets E) \to \prod_\alpha C^M(u_\alpha(U_\alpha),V)
                \\
                s\mapsto \on{pr}_2\circ \psi_\alpha \circ s \circ u_\alpha^{-1}
        \end{gather*}
        where $B\supseteq U_\alpha \East{u_\alpha}{} u_\alpha(U_\alpha)\subset W$ is a
        $C^M$-atlas for $B$ which we assume to be modelled on a convenient vector
        space $W$, and where
        $\psi_\alpha :E|_{U_\alpha}\to U_\alpha\times V$ form a vector bundle atlas over
        charts $U_\alpha$ of $B$. 
        
        \begin{lemma*}
                For a $C^M$ vector bundle $E\to B$ a curve $c:\mathbb R\to C^M(B\gets E)$ is
                $C^M$ if and only if $c^\wedge :\mathbb R\times B\to E$ is $C^M$. 
        \end{lemma*}
        
        \begin{demo}{Proof}
                By the description of the structure on $C^M(B\gets E)$ we may assume that
                $B$ is $c^\infty$-open in a convenient vector space $W$ and that $E=B\times V$. 
                Then $C^M(B\gets B\times V)\cong C^M(B,V)$. Then the statement follows from the
                exponential law \thetag{\ref{nmb:5.3}}. 
                \qed\end{demo}
        
        An immediate consequence is the following: If $U\subset E$ is an open
        neighborhood of $s(B)$ for a section $s$, $F\to B$ is another vector bundle and if $f:U\to
        F$ is a fiber respecting $C^M$ mapping, then $f_*:C^M(B\gets U)\to
        C^M(B\gets F)$ is $C^M$ on the open neighborhood $C^M(B\gets U)$ of $s$ in
        $C^M(B\gets E)$. We have $(d(f_*)(s)v)_x = d(f|_{U\cap E_x})(s(x))(v(x))$.

        \begin{theorem}\label{nmb:6.3}
                Let $M=(M_k)$ be a non-quasianalytic DC-weight sequence of moderate growth. 
                Let $A$ and $B$ be finite dimensional $C^M$ manifolds with $A$ compact. 
                Then the space $C^M(A,B)$ of all $C^M$-mappings $A\to B$ is a
                $C^M$-manifold modelled on convenient vector spaces $C^M(A\gets f^*TB)$ of
                $C^M$ sections of pullback bundles along $f:A\to B$. 
                Moreover, a curve $c:\mathbb R\to C^M(A,B)$ is $C^M$ if and only if 
                $c^\wedge :\mathbb R\times A\to B$ is $C^M$. 
        \end{theorem}
        
        \begin{demo}{Proof}
                Choose a $C^M$ Riemannian metric on $B$ which exists since we have $C^M$
                partitions of unity. $C^M$-vector fields have $C^M$-flows by
                \cite{Komatsu80}; applying this to the geodesic spray 
                we get the $C^M$ exponential mapping 
                $\exp: TB\supseteq U\to B$ of this Riemannian metric,
                defined on a suitable open neighborhood of the zero section. 
                We may assume that $U$ is chosen in such a way that
                $(\pi_B,\exp):U\to B\times B$ is a $C^M$ diffeomorphism onto
                an open neighborhood $V$ of the diagonal, 
                by the $C^M$ inverse function theorem
                due to \cite{Komatsu79}. 
                
                For $f\in C^M(A,B)$ we consider the pullback vector bundle
                \[\xymatrix{
                        A\times_BTB \ar@{=}[r] & f^*TB \ar[r]^{\pi_B^*f} \ar[d]_{f^*\pi_B} & 
                        TB \ar[d]^{\pi_B} \\
                        & A \ar[r]^f & B
                }\]
                Then $C^M(A\gets f^*TB)$ is canonically isomorphic to the space
                $C^M(A,TB)_f:= \{h\in C^M(A,TB):\pi_B\circ h=f\}$ 
                via $s\mapsto (\pi_B^*f)\circ s$ and $(\on{Id}_A,h)\DOTSB\kern.2em\setbox0=\hbox{$\leftarrow$\kern-.15em\raise0.1ex\hbox{$\shortmid$}}\box0\kern.3em h$. 
                Now let 
                \begin{gather*} 
                        U_f :=\{g\in C^M(A,B):(f(x),\; g(x))\in V
                        \text{ for all }x\in A\},\\
                        u_f:U_f\to C^M(A\gets f^*TB),\\
                        u_f(g)(x) = (x,\exp_{f(x)}^{-1}(g(x))) = (x,((\pi_B,\exp)^{-1}\circ(f,g))(x)).
                \end{gather*}
                Then $u_f$ is a bijective mapping from $U_f$ onto the set
                $\{s\in C^M(A\gets f^*TB): s(A)\subseteq f^*U=(\pi_B^*f)^{-1}(U)\}$, whose 
                inverse is 
                given by $u_f^{-1}(s) = \exp\circ(\pi_B^*f)\circ s$, where we view 
                $U \to B$ as fiber bundle. The push forward $u_f$ is $C^M$ since it maps
                $C^M$-curves to $C^M$-curves by lemma \thetag{\ref{nmb:6.2}}. 
                The set $u_f(U_f)$ is open in
                $C^M(A\gets f^*TB)$ for the topology described above in \thetag{\ref{nmb:6.2}}.
                
                Now we consider the atlas $(U_f,u_f)_{f\in C^M(A,B)}$ for
                $C^M(A,B)$. Its chart change mappings are given 
                for $s\in u_g(U_f\cap U_g)\subseteq C^M(A\gets g^*TB)$ by 
                \begin{align*} 
                        (u_f\circ u_g^{-1})(s) &= (\on{Id}_A,(\pi_B,\exp)^{-1}\circ(f,\exp\circ(\pi_B^*g)\circ s)) \\
                        &= (\tau_f^{-1}\circ\tau_g)_*(s),
                \end{align*}
                where $\tau_g(x,Y_{g(x)}) := (x,\exp_{g(x)}(Y_{g(x)}))$
                is a $C^M$ diffeomorphism 
                $\tau_g:g^*TB \supseteq g^*U \to (g\times \on{Id}_B)^{-1}(V)\subseteq A\times B$
                which is fiber respecting over $A$. 
                The chart change $u_f\circ u_g^{-1} = (\tau_f^{-1}\circ \tau_g)_*$ is defined on an open
                subset and it is also $C^M$ 
                since it respects $C^M$-curves.
                
                Finally for the topology on $C^M(A,B)$ 
                we take the identification topology from this atlas (with the 
                $c^\infty$-topologies on the modeling spaces), 
                which is obviously finer than the
                compact-open topology and thus Hausdorff.
                
                The equation $u_f\circ u_g^{-1} = (\tau_f^{-1}\circ \tau_g)_*$ shows that
                the $C^M$ structure does not depend on the choice of the
                $C^M$ Riemannian metric on $B$.
                
                The statement on $C^M$-curves follows from lemma \thetag{\ref{nmb:6.2}}. 
                \qed\end{demo}
        
        \begin{corollary}\label{nmb:6.4}
                Let $A_1,A_2$ and $B$ be finite dimensional $C^M$ manifolds with $A_1$ and
                $A_2$ compact. 
                Then composition 
                \[
                C^M(A_2,B) \times C^M(A_1,A_2) \to C^M(A_1,B), \quad (f,g) \mapsto f\circ g
                \]
                is $C^M$. 
                However, if $N=(N_k)$ is another non-quasianalytic DC-weight sequence of
                moderate growth with 
                $(N_k/M_k)^{1/k} \searrow 0$ 
                then composition is {\bf not} $C^N$.
        \end{corollary}
        
        \demo{Proof}
        Composition maps $C^M$-curves to $C^M$-curves, so it is $C^M$. 
        
        Let $A_1=A_2=S^1$ and $B=\mathbb R$. Then by \eqref{incl} there exists $f\in
        C^M(S^1,\mathbb R)\setminus C^N(S^1,\mathbb R)$. 
        We consider $f:\mathbb R\to \mathbb R$ periodic.
        The universal covering space of $C^M(S^1,S^1)$ consists of all 
        $2\pi\mathbb Z$-equivariant mappings in $C^M(\mathbb R,\mathbb R)$, 
        namely the space of all $g+\on{Id}_{\mathbb R}$ for $2\pi$-periodic 
        $g\in C^M$. Thus $C^M(S^1,S^1)$ is a real analytic manifold and 
        $t\mapsto (x\mapsto x+t)$ induces a real analytic curve $c$ in $C^M(S^1,S^1)$. 
        But $f_*\circ c$ is not $C^N$ since:
        \begin{align*}
                \frac{(\partial_t^k|_{t=0}(f_*\circ c)(t))(x)}{k!\rho^k N_k} =
                \frac{\partial_t^k|_{t=0} f(x+t)}{k!\rho^k N_k} = \frac{f^{(k)}(x)}{k!\rho^k N_k}
        \end{align*}
        which is unbounded for $x$ in a suitable compact set and for all $\rho>0$ since $f\notin C^N$.
        \qed\enddemo
        
        \begin{theorem}\label{nmb:6.5}
                Let $M=(M_k)$ be a non-quasianalytic DC-weight sequence of moderate growth. 
                Let $A$ be a compact ($\implies$ finite dimensional) $C^M$ manifold.
                Then the group $\on{Diff}^M(A)$ of all $C^M$-diffeomorphisms of $A$ is an
                open subset of the $C^M$ manifold $C^M(A,A)$. Moreover, it is
                a $C^M$-regular $C^M$ Lie group: 
                Inversion and composition are $C^M$. Its Lie algebra
                consists of all $C^M$-vector fields on $A$, with the negative of the usual
                bracket as Lie bracket. 
                The exponential mapping is $C^M$. It is not surjective onto any neighborhood
                of $\on{Id}_A$. 
        \end{theorem}
        
        Following \cite{KM97r}, see also \cite[38.4]{KM97}, 
        a $C^M$-Lie group $G$ with Lie algebra $\mathfrak g=T_eG$ 
        is called $C^M$-regular if the following holds:
        \begin{itemize}
        \item 
                For each $C^M$-curve 
                $X\in C^M(\mathbb R,\mathfrak g)$ there exists a $C^M$-curve 
                $g\in C^M(\mathbb R,G)$ whose right logarithmic derivative is $X$, i.e.,
                \[
                \begin{cases} g(0) &= e \\
                        \partial_t g(t) &= T_e(\mu^{g(t)})X(t) = X(t).g(t)
                \end{cases} 
                \]
                The curve $g$ is uniquely determined by its initial value $g(0)$, if it
                exists.
        \item
                Put $\on{evol}^r_G(X)=g(1)$ where $g$ is the unique solution required above. 
                Then $\on{evol}^r_G: C^M(\mathbb R,\mathfrak g)\to G$ is required to be
                $C^M$ also. 
        \end{itemize}
        
        \demo{Proof}
        The group $\on{Diff}^M(A)$ is open in $C^M(A,A)$ since it is open in the
        coarser $C^1$ compact open topology, see \cite[43.1]{KM97}. 
        So $\on{Diff}^M(A)$ is a $C^M$-manifold and composition is $C^M$ by
        \thetag{\ref{nmb:6.3}} and \thetag{\ref{nmb:6.4}}. To show that inversion is $C^M$ let $c$ be a
        $C^M$-curve in $\on{Diff}^M(A)$. By \thetag{\ref{nmb:6.3}} the map 
        $c^\wedge: \mathbb R\times A\to A$ is
        $C^M$ and $(\on{inv}\circ c)^\wedge:\mathbb R\times A\to A$ satisfies the finite
        dimensional implicit equation 
        $c^\wedge(t,(\on{inv}\circ c)^\wedge(t,x))=x$ for $x\in A$. By the finite
        dimensional $C^M$ implicit function theorem \cite{Komatsu79} the mapping 
        $(\on{inv}\circ c)^\wedge$ is locally $C^M$ and thus $C^M$. By \thetag{\ref{nmb:6.3}}
        again, $\on{inv}\circ c$ is a $C^M$-curve in $\on{Diff}^M(A)$.
        So $\on{inv}:\on{Diff}^M(A)\to \on{Diff}^M(A)$ is $C^M$. 
        The Lie algebra of $\on{Diff}^M(A)$ is the convenient vector space of all 
        $C^M$-vector fields on $A$, with the negative of the usual Lie bracket
        (compare with the proof of \cite[43.1]{KM97}). 
        
        To show that $\on{Diff}^M(A)$ is a $C^M$-regular Lie group, we choose a
        $C^M$-curve in the space of $C^M$-curves in the Lie algebra of all $C^M$
        vector fields on $A$, $c:\mathbb R\to C^M(\mathbb R,C^M(A\gets TA))$. By
        lemma \thetag{\ref{nmb:6.2}} $c$ corresponds to a $\mathbb R^2$-time-dependent $C^M$
        vector field $c^{\wedge \wedge }:\mathbb R^2\times A\to TA$. Since $C^M$-vector
        fields have $C^M$-flows and since $A$ is compact, 
        $\on{evol}^r(c^\wedge (s))(t) = \on{Fl}^{c^\wedge (s)}_t$ is $C^M$ in all
        variables by \cite{Komatsu80}. 
        Thus $\on{Diff}^M(A)$ is a $C^M$-regular $C^M$ Lie group.
        
        The exponential mapping is $\on{evol}^r$ applied to constant curves in the
        Lie algebra, i.e., it consists of flows of autonomous $C^M$ vector fields. 
        That the exponential map is not surjective onto any $C^M$-neighborhood of
        the identity follows from \cite[43.5]{KM97} for $A=S^1$. This example can
        be embedded into any compact manifold, see \cite{Grabowski88}. 
        \qed\enddemo

        \section{\label{nmb:7} Appendix. Calculus beyond Banach spaces} 
        The traditional differential calculus works 
        well for finite dimensional vector spaces and for Banach spaces. For 
        more general locally convex spaces 
        we sketch here the convenient approach as explained in 
        \cite{FK88} and \cite{KM97}.
        The main difficulty is that composition of 
        linear mappings stops to be jointly continuous at the level of Banach 
        spaces, for any compatible topology. 
        We use the notation of \cite{KM97} and this is the
        main reference for the whole appendix. We list results in the order in
        which one can prove them, without proofs for which we refer to \cite{KM97}. 
        This should explain how to use these results.
        
        \subsection{\label{nmb:7.1}The $c^\infty$-topology}
        Let $E$ be a 
        locally convex vector space. A curve $c:\mathbb R\to E$ is called 
        {\it smooth} or $C^\infty$ if all derivatives exist and are 
        continuous - this is a concept without problems. Let 
        $C^\infty(\mathbb R,E)$ be the space of smooth functions. It can be 
        shown that the set $C^\infty(\mathbb R,E)$ does not depend on the locally convex 
        topology of $E$, only on its associated bornology (system of bounded 
        sets).
        
        The final topologies with respect to the following sets of mappings into E coincide:
        \begin{enumerate}
        \item $C^\infty(\mathbb R,E)$.
        \item The set of all Lipschitz curves 
                (so that $\{\frac{c(t)-c(s)}{t-s}:t\neq s\}$ 
                is bounded in $E$). 
        \item The set of injections $E_B\to E$ where $B$ runs through all bounded 
                absolutely convex subsets in $E$, and where 
                $E_B$ is the linear span of $B$ equipped with the Minkowski 
                functional $\|x\|_B:= \inf\{\lambda>0:x\in\lambda B\}$.
        \item The set of all Mackey-convergent sequences $x_n\to x$ 
                (there exists a sequence 
                $0<\lambda_n\nearrow\infty$ with $\lambda_n(x_n-x)$ bounded).
        \end{enumerate}
        This topology is called the $c^\infty$-topology on $E$ and we write 
        $c^\infty E$ for the resulting topological space. In general 
        (on the space $\mathcal{D}$ of test functions for example) it is finer 
        than the given locally convex topology, it is not a vector space 
        topology, since scalar multiplication is no longer jointly 
        continuous. The finest among all locally convex topologies on $E$ 
        which are coarser than $c^\infty E$ is the bornologification of the 
        given locally convex topology. If $E$ is a Fr\'echet space, then 
        $c^\infty E = E$. 
        
        \subsection{\label{nmb:7.2}Convenient vector spaces} 
        A locally convex vector space 
        $E$ is said to be a {\it convenient 
                vector space} if one of the following equivalent
        conditions is satisfied (called $c^\infty$-completeness):
        \begin{enumerate}
        \item For any $c\in C^\infty(\mathbb R,E)$ the (Riemann-) integral 
                $\int_0^1c(t)dt$ exists in $E$.
        \item Any Lipschitz curve in $E$ is locally Riemann integrable.
        \item A curve $c:\mathbb R\to E$ is smooth if and only if $\lambda\circ c$ is 
                smooth for all $\lambda\in E^*$, where $E^*$ is the dual consisting 
                of all continuous linear functionals on $E$. Equivalently, 
                we may use the dual $E'$ consisting of 
                all bounded linear functionals.
        \item Any Mackey-Cauchy-sequence (i.\ e.\  $t_{nm}(x_n-x_m)\to 0$ 
                for some $t_{nm}\to \infty$ in $\mathbb R$) converges in $E$. 
                This is visibly a mild completeness requirement.
        \item If $B$ is bounded closed absolutely convex, then $E_B$ 
                is a Banach space.
        \item If $f:\mathbb R\to E$ is scalarwise $\operatorname{\mathcal Lip}^k$, then $f$ is 
                $\operatorname{\mathcal Lip}^k$, for $k>1$.
        \item If $f:\mathbb R\to E$ is scalarwise $C^\infty$ then $f$ is 
                differentiable at 0.
        \item If $f:\mathbb R\to E$ is scalarwise $C^\infty$ then $f$ is 
                $C^\infty$.
        \end{enumerate}
        Here a mapping $f:\mathbb R\to E$ is called $\operatorname{\mathcal Lip}^k$ if all 
        derivatives up to order $k$ exist and are Lipschitz, locally on 
        $\mathbb R$. That $f$ is scalarwise $C^\infty$ means $\lambda\circ f$ is $C^\infty$ 
        for all continuous linear functionals on $E$.
        
        \subsection{\label{nmb:7.3}Smooth mappings} 
        Let $E$, $F$, and $G$ be convenient vector spaces, 
        and let $U\subset E$ be $c^\infty$-open. 
        A mapping $f:U\to F$ is called {\it smooth} or 
        $C^\infty$, if $f\circ c\in C^\infty(\mathbb R,F)$ for all 
        $c\in C^\infty(\mathbb R,U)$.
        {\it
                The main properties of smooth calculus are the following.
                \begin{enumerate}
                \item For mappings on Fr\'echet spaces this notion of smoothness 
                        coincides with all other reasonable definitions. Even on 
                        $\mathbb R^2$ this is non-trivial.
                \item Multilinear mappings are smooth if and only if they are 
                        bounded.
                \item If $f:E\supseteq U\to F$ is smooth then the derivative 
                        $df:U\times E\to F$ is 
                        smooth, and also $df:U\to L(E,F)$ is smooth where $L(E,F)$ 
                        denotes the space of all bounded linear mappings with the 
                        topology of uniform convergence on bounded subsets.
                \item The chain rule holds.
                \item The space $C^\infty(U,F)$ is again a convenient vector space 
                        where the structure is given by the obvious injection
                        \[
                        C^\infty(U,F) \East{C^\infty(c,\ell)}{} 
                        \negthickspace\negthickspace\negthickspace\negthickspace\negthickspace
                        \negthickspace\negthickspace\negthickspace
                        \prod_{c\in C^\infty(\mathbb R,U), \ell\in F^*} 
                        \negthickspace\negthickspace\negthickspace\negthickspace\negthickspace
                        \negthickspace\negthickspace\negthickspace
                        C^\infty(\mathbb R,\mathbb R),
                        \quad f\mapsto (\ell\circ f\circ c)_{c,\ell},
                        \]
                        where $C^\infty(\mathbb R,\mathbb R)$ carries the topology of compact 
                        convergence in each derivative separately.
                \item The exponential law holds: For $c^\infty$-open $V\subset F$, 
                        \[
                        C^\infty(U,C^\infty(V,G)) \cong C^\infty(U\times V, G)
                        \]
                        is a linear diffeomorphism of convenient vector spaces. Note 
                        that this is the main assumption of variational calculus.
                \item A linear mapping $f:E\to C^\infty(V,G)$ is smooth (bounded) if 
                        and only if $E \East{f}{} C^\infty(V,G) \East{\on{ev}_v}{} G$ is smooth 
                        for each $v\in V$. This is called the smooth uniform 
                        boundedness theorem \cite[5.26]{KM97}.
                \item The following canonical mappings are smooth.
                        \begin{align*}
                                &\operatorname{ev}: C^\infty(E,F)\times E\to F,\quad 
                                \operatorname{ev}(f,x) = f(x)\\
                                &\operatorname{ins}: E\to C^\infty(F,E\times F),\quad
                                \operatorname{ins}(x)(y) = (x,y)\\
                                &(\quad)^\wedge :C^\infty(E,C^\infty(F,G))\to C^\infty(E\times F,G)\\
                                &(\quad)^\vee :C^\infty(E\times F,G)\to C^\infty(E,C^\infty(F,G))\\
                                &\operatorname{comp}:C^\infty(F,G)\times C^\infty(E,F)\to C^\infty(E,G)\\
                                &C^\infty(\quad,\quad):C^\infty(F,F_1)\times C^\infty(E_1,E)\to 
                                C^\infty(C^\infty(E,F),C^\infty(E_1,F_1))\\
                                &\qquad (f,g)\mapsto(h\mapsto f\circ h\circ g)\\
                                &\prod:\prod C^\infty(E_i,F_i)\to C^\infty(\prod E_i,\prod F_i)
                        \end{align*}
                \end{enumerate}
        }
        
        \subsection{\label{nmb:7.4}Remarks} 
        Note that the conclusion of 
        \thetag{\ref{nmb:7.3}.6} is the starting point of the classical calculus of 
        variations, where a smooth curve in a space of functions was assumed 
        to be just a smooth function in one variable more. It is also the source of the name convenient 
        calculus.
        This and some other obvious properties already determines the convenient calculus.
        
        There are, however, smooth mappings which are not continuous. This is 
        unavoidable and not so horrible as it might appear at first sight. 
        For example the evaluation $E\times E^*\to\mathbb R$ is jointly continuous if 
        and only if $E$ is normable, but it is always smooth. Clearly smooth 
        mappings are continuous for the $c^\infty$-topology.

        \section{\label{nmb:8} Calculus of holomorphic mappings}
        
        \subsection{\label{nmb:8.1}Holomorphic curves} 
        Let $E$ be a complex locally convex vector 
        space whose underlying real space is convenient -- this will be 
        called convenient in the sequel. Let $\mathbb D\subset \mathbb C$ be the 
        open unit disk and let us denote by $\mathcal{H}(\mathbb D,E)$ the space of 
        all mappings $c:\mathbb D\to E$ such that $\lambda\circ c:\mathbb D\to \mathbb C$ is 
        holomorphic for each continuous complex-linear functional $\lambda$ on 
        $E$. Its elements will be called the holomorphic curves.
        
        If $E$ and $F$ are convenient complex vector spaces (or 
        $c^\infty$-open sets therein), a mapping 
        $f:E\to F$ is called {\it holomorphic} if $f\circ c$ is a holomorphic 
        curve in $F$ for each holomorphic curve $c$ in $E$. Obviously $f$ is 
        holomorphic if and only if $\lambda\circ f:E\to \mathbb C$ is holomorphic for 
        each complex linear continuous (equivalently: bounded) functional $\lambda$ on $F$. Let 
        $\mathcal{H}(E,F)$ denote the space of all holomorphic mappings from $E$ to 
        $F$. 
        
        \begin{lemma}\label{nmb:8.2}(Hartog's theorem) 
                Let $E_k$ for 
                $k=1,2$ and $F$ be complex convenient vector spaces and let 
                $U_k\subset E_k$ be $c^\infty$-open. A mapping $f:U_1\times U_2\to F$ is 
                holomorphic if and only if it is separately holomorphic (i.\ e.\ 
                $f(\quad,y)$ and $f(x,\quad)$ are holomorphic for all $x\in U_1$ and 
                $y\in U_2$).
        \end{lemma}
        
        This implies also that in finite dimensions we have recovered the 
        usual definition.
        
        \begin{lemma}\label{nmb:8.3}
                If $f:E\supset U\to F$ is holomorphic 
                then $df:U\times E\to F$ exists, is holomorphic and $\mathbb C$-linear in 
                the second variable. 
                
                A multilinear mapping is holomorphic if and only if it is bounded.
        \end{lemma}
        
        \begin{lemma}\label{nmb:8.4}
                If $E$ and $F$ are Banach spaces and $U$ 
                is open in $E$, then for a mapping $f:U\to F$ the following 
                conditions are equivalent:
                \begin{enumerate}
                \item $f$ is holomorphic.
                \item $f$ is locally a convergent series of homogeneous continuous 
                        polynomials.
                \item $f$ is $\mathbb C$-differentiable in the sense of Fr\'echet.
                \end{enumerate} 
        \end{lemma}
        
        \begin{lemma}\label{nmb:8.5}
                Let $E$ and $F$ be convenient vector 
                spaces. A mapping $f:E\to F$ is holomorphic if and only if it is 
                smooth and its derivative in each point is $\mathbb C$-linear.
        \end{lemma}
        
        An immediate consequence of this result is that $\mathcal{H}(E,F)$ is a 
        closed linear subspace of $C^\infty(E_{\mathbb R},F_{\mathbb R})$ and so it 
        is a convenient vector space if $F$ is one, by \thetag{\ref{nmb:7.3}.5}. 
        The chain rule follows from \thetag{\ref{nmb:7.3}.4}. 
        
        \begin{theorem}\label{nmb:8.6}
                The category of convenient complex 
                vector spaces and holomorphic mappings between them is cartesian 
                closed, i.\ e.
                \begin{displaymath}
                        \mathcal{H}(E\times F,G) \cong \mathcal{H}(E,\mathcal{H}(F,G)).\end{displaymath}
        \end{theorem}
        
        An immediate consequence of this is again that all canonical 
        structural mappings as in \thetag{\ref{nmb:7.3}.8} are holomorphic.
        
        \section{\label{nmb:9} Calculus of real analytic mappings}
        
        \subsection{\label{nmb:9.1}} We now sketch the cartesian closed 
        setting to real analytic mappings in infinite dimension following the 
        lines of the Fr\"olicher--Kriegl calculus, as it is presented in 
        \cite{KM97}. Surprisingly enough one has to deviate 
        from the most obvious notion of real analytic curves in order to get 
        a meaningful theory, but again convenient vector spaces turn out to 
        be the right kind of spaces.
        
        \subsection{\label{nmb:9.2} Real analytic curves} Let $E$ be a real 
        convenient vector space with continuous dual $E^*$. A curve $c:\mathbb R\to E$ is 
        called {\it real analytic} if $\lambda\circ c:\mathbb R\to \mathbb R$ is real 
        analytic for each $\lambda\in E^*$. 
        It turns out that the set of these curves depends only on the 
        bornology of $E$.
        Thus we may use the dual $E'$ consisting of all bounded 
        linear functionals in the definition.

        In contrast a curve is called {\it strongly real analytic} if it 
        is locally given by power series which converge in the topology of 
        $E$. They can be extended to germs of holomorphic curves along $\mathbb R$ 
        in the complexification $E_{\mathbb C}$ of $E$. If the dual $E^*$ of $E$ 
        admits a Baire topology which is compatible with the duality, then 
        each real analytic curve in $E$ is in fact topologically real analytic 
        for the bornological topology on $E$.
        
        \subsection{\label{nmb:9.3} Real analytic mappings} Let $E$ and $F$ be 
        convenient vector spaces. Let $U$ be a $c^\infty$-open set in $E$. A 
        mapping $f:U\to F$ is called {\it real analytic} if and only if it is 
        smooth (maps smooth curves to smooth curves) and maps real analytic 
        curves to real analytic curves. 
        
        Let $C^\omega(U,F)$ denote the space of all real analytic mappings. 
        We equip the space
        $C^\omega(U,\mathbb R)$ of all real analytic functions 
        with the initial topology
        with respect to the families of mappings
        \begin{gather*} C^\omega(U,\mathbb R) \East{{c^*}}{}C^\omega(\mathbb R,\mathbb R),\text{
                        for all }c\in C^\omega(\mathbb R,U)\\ 
                C^\omega(U,\mathbb R) \East{{c^*}}{}C^\infty(\mathbb R,\mathbb R),\text{
                        for all }c\in C^\infty(\mathbb R,U),
        \end{gather*}
        where $C^\infty(\mathbb R,\mathbb R)$ carries the topology of compact 
        convergence in each derivative separately, and 
        where $C^\omega(\mathbb R,\mathbb R)$ is equipped with the final locally 
        convex topology 
        with respect to the embeddings (restriction mappings) of all spaces 
        of holomorphic mappings from a neighborhood $V$ of $\mathbb R$ in 
        $\mathbb C$ mapping $\mathbb R$ to $\mathbb R$, and each of these spaces 
        carries the topology of compact convergence.
        
        Furthermore we equip the space
        $C^\omega(U,F)$ with the initial topology with respect to
        the family of mappings 
        \begin{displaymath}
                C^\omega(U,F) \East{{\lambda_*}}{}C^\omega(U,\mathbb R),\text{ for all }\lambda\in
                F^*.
        \end{displaymath}
        It turns out that this is again a convenient space.
        
        \begin{theorem}\label{nmb:9.4}In the setting of \thetag{\ref{nmb:9.3}} a mapping 
                $f:U\to F$ is real analytic if and only if it is smooth and is real 
                analytic along each affine line in $E$.
        \end{theorem}
        
        \begin{lemma}\label{nmb:9.5}The space $L(E,F)$ of all bounded linear 
                mappings is a closed linear subspace of $C^\omega(E,F)$. A mapping 
                $f:U\to L(E,F)$ is real analytic if and only if $\on{ev}_x\circ f:U\to F$ is 
                real analytic for each point $x\in E$.
        \end{lemma}
        
        \begin{theorem}\label{nmb:9.6}The category of convenient spaces and 
                real analytic mappings is cartesian closed. So the equation
                \begin{displaymath}
                        C^\omega(U,C^\omega(V,F))\cong C^\omega(U\times V,F)\end{displaymath}
                is valid for all $c^\infty$-open sets $U$ in $E$ and $V$ in $F$, 
                where $E$, $F$, and $G$ are convenient vector spaces.
        \end{theorem}
        
        This implies again that all structure mappings as in \thetag{\ref{nmb:7.3}.8} are 
        real analytic. Furthermore the differential operator 
        \begin{displaymath}
                d:C^\omega(U,F)\to C^\omega(U,L(E,F))
        \end{displaymath}
        exists, is unique and real 
        analytic. Multilinear mappings are real analytic if and only if they 
        are bounded. 
        
        \begin{theorem}[Real analytic uniform boundedness principle]\label{nmb:9.7}
                A linear mapping $f:E\to C^\omega(V,G)$ is real analytic (bounded) if 
                and only if $E \East{f}{} C^\omega(V,G) \East{\on{ev}_v}{} G$ is real
                analytic (bounded) for each $v\in V$. 
        \end{theorem}

\def\cprime{$'$}
\providecommand{\bysame}{\leavevmode\hbox to3em{\hrulefill}\thinspace}
\providecommand{\MR}{\relax\ifhmode\unskip\space\fi MR }
\providecommand{\MRhref}[2]{%
  \href{http://www.ams.org/mathscinet-getitem?mr=#1}{#2}
}
\providecommand{\href}[2]{#2}

\end{document}